\documentclass[leqno,11pt]{article}
\usepackage{latexsym}
\usepackage{amssymb}
\usepackage{amsfonts}
\usepackage{amsmath}
\usepackage{epsfig}
\usepackage{color}
\usepackage{pifont}

\allowdisplaybreaks
\makeatletter
\long\def\unmarkedfootnote#1{{\long\def\@makefntext##1{##1}\footnotetext{#1}}}
\makeatother

\newtheorem{definition}{Definition}[section]

\newtheorem{lemma}[definition]{Lemma}

\newtheorem{theorem}[definition]{Theorem}

\newtheorem{proposition}[definition]{Proposition}

\newtheorem{corollary}[definition]{Corollary}

\newtheorem{remark}[definition]{Remark}

\def\o{\Omega}

\def\m2{|\Omega | /2}
\def\M2{\frac{|\Omega |}{2}}
\def\u+{u_+^*}

\def\-p{\overline{p}}

\def\w0{{W_0^{1,p}(\Omega)}}

\def\R{\mathbb R}

\def\ep{\varepsilon}
\def\rN{\mathbb R^N}
\def\bu{{\bf u}}

\newcommand{\bm}[1]{{\boldsymbol{\rm #1}}}
\def\rn{{{\R}^n}}
\def\rN{{{\R}^N}}
\def\rNn{{{\R}^{Nn}}}

\newcommand{\hh}{{\cal H}^{n-1}}
\newcommand{\medint}{-\kern  -,395cm\int}
\newcommand{\medintinrigo}{-\kern  -,315cm\int}
\newcommand{\medelle}{-\kern  -,235cm L}
\newcommand{\medellenrigo}{-\kern  -,180cm L}
\newcommand{\qed}{\thinspace\null\nobreak\hfill
\hbox{\vbox{\kern-.2pt\hrule height.2pt
depth.2pt\kern-.2pt\kern-.2pt \hbox to1.8mm {\kern-.2pt\vrule
width.4pt \kern-.2pt\raise1.8mm\vbox to.2pt{} \lower0pt\vtop
to.2pt{}\hfil\kern-.2pt \vrule
width.4pt\kern-.2pt}\kern-.2pt\kern-.2pt \hrule height.2pt
depth.2pt \kern-.2pt}}\par\medbreak}

\title{Optimal second-order regularity  for the $p$-Laplace system
} 
\numberwithin{equation}{section}

\author{
  Andrea Cianchi\\
 {\it Dipartimento di Matematica e Informatica \lq\lq U. Dini", Universit\`a di Firenze}\\ {\it Viale Morgagni 67/A, 50134 Firenze, Italy} \\{\it  e-mail: cianchi@unifi.it}
\bigskip
\\
  Vladimir G. Maz'ya \\
  {\it   Department of Mathematics, Link\"oping University, SE-581
83 Link\"oping, Sweden}
  \\ and \\
{\it  RUDN University}\\ {\it
6 Miklukho-Maklay St, Moscow, 117198, Russia}
\\ {\it e-mail: vladimir.mazya@liu.se}
}

\date{}

\begin{document}
\maketitle

\begin{abstract}
Second-order estimates are established for solutions to the $p$-Laplace system with right-hand side in $L^2$. The nonlinear expression of the gradient under the divergence operator is shown to belong to $W^{1,2}$, and hence to enjoy the best possible degree of regularity. Moreover, its norm in $W^{1,2}$ is proved to be equivalent to the norm of the right-hand side in $L^2$. Our global results apply to solutions to both  Dirichlet and Neumann problems, and entail minimal regularity of the boundary of the domain. 
In particular, our conclusions hold for arbitrary bounded convex domains.
Local estimates for local solutions are provided as well.
\end{abstract}

\unmarkedfootnote {

%
%
}

\section{Introduction}\label{intro}

We are concerned with second-order differentiability properties of solutions to 
 the $p$-Laplace system
\begin{equation}\label{localeq}
- {\rm {\bf  div}} (|\nabla \bu|^{p-2}\nabla {\bf u} ) = {\bf f} \quad {\rm in}\,\,\, \o\,.
\end{equation}
Here, $\Omega$ is an open set in $\rn$,   $n \geq 2$,  the function ${\bf f} : \Omega \to \mathbb R^N$, $N \geq 1$, is  given, and  $\bu: \o \to \rN$ is the unknown.
\par Local solutions to system \eqref{localeq} are included in our analysis. However, our main focus is on Dirichlet boundary value problems, obtained on coupling system \eqref{localeq} with the boundary condition 
\begin{equation}\label{dircond}
\bu =0 \qquad \hbox{on $\partial \o$,}
\end{equation}
and on Neumann problems, corresponding to the boundary condition
\begin{equation}\label{neucond}
\frac{\partial \bu}{\partial \bm \nu} =0 \qquad \hbox{on $\partial \o$,}
\end{equation}
where $\bm \nu$ denotes the outward unit vector on $\partial \Omega$. Obviously, the compatibility assumption
\begin{equation}\label{0mean}
\int _\o {\bf f}\, dx =0
\end{equation}
 has also to be imposed  on $\bm f$ when  \eqref{neucond} is in force. 
%
%
\par
It is clear that, if $\bm u$ is a solution to system \eqref{localeq} such that $|\nabla \bu|^{p-2}\nabla {\bf u} \in W^{1,2}(\o, \rNn)$, then $\bm f \in L^2(\o, \rN)$. Our main result amounts to a converse of this fact, and furnishes a nonlinear counterpart of  the classical $L^2$-coercivity theory  for linear problems, that can be traced back to \cite{Bernstein} ($n=2$) and \cite{Schauder} (arbitrary $n$).
 Namely, it asserts that if $\bm f \in L^2(\o, \rN)$, then $|\nabla \bu|^{p-2}\nabla {\bf u} \in W^{1,2}(\o, \rNn)$, and, moreover, the norms $\|\bm f \|_{L^2(\o, \rN)}$ and  $\||\nabla \bu|^{p-2}\nabla {\bf u}\|_{W^{1,2}(\o, \rNn)}$ are equivalent.
Of course, this loose statement requires some specifications. 
\par To begin with, the 
equivalence of the relevant global norms  for 
%
%
solutions to either the Dirichlet problem \eqref{localeq}$+$\eqref{dircond}, or the Neumann problem \eqref{localeq}$+$\eqref{neucond}, can only hold if  $\partial \o$ is sufficiently regular. In this connection, let us stress that, athough our results are new even for smooth domains, we strive for weakest possible regularity assumptions on $\partial \o$. As will be demonstrated by apropos examples, the hypotheses to be made on $\partial \o$ are virtually minimal. On the other hand, 
 local solutions obviously just admit local estimates for the  $W^{1,2}$-norm of $|\nabla \bu|^{p-2}\nabla {\bf u}$, which do not entail any assumption on the domain $\o$. However, besides a local $L^2$-norm of $\bm f$, they involve an additional lower-order term depending on $\nabla \bm u$. 
\par
The very definition of solution to system \eqref{localeq} calls for  an elucidation.  Indeed, membership of $\bm f$ in $L^2(\o, \rN)$ does not ensure that weak solutions be well defined, if $p$ is too small -- precisely, if $p< \frac{2n}{n+2}$, the H\"older conjugate of the Sobolev conjugate of the exponent $2$. As a consequence, an even weaker notion of solution has to be employed. In the scalar case, namely  when $N=1$ and hence system \eqref{localeq} reduces to a single equation, diverse (a posteriori equivalent) definitions of solutions are available in the literature that allow for right-hand sides $\bm f$ that are merely integrable in $\o$ -- see \cite{BBGGPV, BG, DallA, DM, LM, Ma69, Mu}. The existence of a solution of any of these kinds can be  established, both under the Dirichlet condition 
\eqref{dircond} and the Neumann  condition \eqref{neucond},
  via auxiiary approximating problems, whose  right-hand sides are smooth. The solution turns out to be the pointwise limit of the solutions to the axiliary problems,  and  need not even be weakly differentiable if $p$ is too small, although it is always endowed with a surrogate of a gradient, given by the pointwise limit of the gradients of the approximating solutions. In case of equations, this approximation method applies for every $p>1$.
\par
By contrast, in the case of systems with right-hand sides affected by weak integrability properties, the existence of solutions to boundary value problems obtained as the limit of smooth solutions to approximating problems 
is only known  under the assumption that $p>2-\frac 1n$ -- see  \cite{DHM1, DHM3}. Solutions in this sense are  called approximable solutions hereafter, and are precisely defined in Sections \ref{proofsloc} and \ref{proofs}.  Our results will be established  for the even broader range of exponents $p > \tfrac 32$, a threshold  that  arises in a  crucial differential inequality entering  in our proof. 
\par Such a limitation on $p$ is not required when  equations, instead of  systems, are in question. This is shown in our earlier contribution \cite{cmsecond}, where second-order regularity results in the same spirit are established for the $p$-Laplace (and a somewhat more general quasilinear) equation, for every $p>1$. A simpler differential inequality,  which holds with no restriction on $p$, comes into play in that case. We emphasize that, besides dealing with systems, another advance of the present paper consists in further relaxing the regularity of $\partial \o$. This calls for a new approximation result for  domains $\o$ in the class under consideration, which is needed in the proof. We have also to resort to some additional step in the   approximation scheme for  solutions, since certain global regularity results in case of smooth data and domains, that are well known for equations, are not available for systems.
\par 
Let us mention that classical results on second-order differentiability properties of solutions to the $p$-Laplacian equation or  system are concerned with the expression $|\nabla u|^{\frac{p-2}2}\nabla u$, instead  of $|\nabla u|^{p-2}\nabla u$. They can be traced back to \cite{Uhl} for $p>2$, and \cite{CDiB} for every $p >1$.  Recent developments are in   \cite{BC, Cel, CGM}. 
The differentiability of the expression $|\nabla u|^{{p-2}}\nabla u$ has also been investigated, but under stronger regularity assumptions on the right-hand side than the natural membership in $L^2(\o)$. Furthermore, the available results in this direction either require smoothness of $\partial \o$ \cite{Dasc}, or just deal with local solutions \cite{AKM, Lou}. 
Fractional-order regularity  of the gradient of solutions to quasilinear  equations  of   $p$-Laplacian type has been studied in \cite{SimonJ}, and   recently in \cite{AKM, Cel1, Mi1}.

%
%
%

\section{Main results}\label{introbis}

Our most general global result  for solutions to Dirichlet and Neumann problems entails  a regularity condition on $\partial \o$ depending, roughly speaking, on a local isocapacitary inequality for 
 the integral of the curvature on $\partial \o$. 
\par Precisely,
we assume that $\Omega$ is  a bounded Lipschitz domain, and that the functions of    $(n-1)$ variables that locally describe the boundary of $\o$ are   twice weakly differentiable, briefly $\partial \o \in W^{2,1}$. These assumptions do not yet ensure global second-order regularity of solutions, still in the linear case -- see \cite{Ma67, Ma73}. They have to be refined  as follows. Denote by 
$\mathcal B$ the weak second fundamental form on $\partial \o$, by
$|\mathcal B|$ its norm, and set 
\begin{equation}\label{defK}
\mathcal K_\o(r) =
\displaystyle 
\sup _{
\begin{tiny}
 \begin{array}{c}{E\subset \partial \o \cap B_r(x)} \\
x\in \partial \o
 \end{array}
  \end{tiny}
} \frac{\int _E |\mathcal B|d\hh}{{\rm cap}_{B_1(x)} (E)}\qquad \hbox{for $r\in (0,1)$}\,.
\end{equation}
Here, $B_r(x)$ stands for the ball centered at $x$, with radius $r$,  the notation ${\rm cap}_{B_1(x)}(E)$ is adopted for the capacity of the set $E$ relative to the ball $B_1(x)$, and $\hh$ is the $(n-1)$-dimensional Hausdorff measure.
Then we require that
\begin{equation}\label{capcond}
\lim _{r\to 0^+} \mathcal K _\o(r) < c
\end{equation}
for a suitable constant $c=c(n,N,p,L_\o, d_\o)$,
where, 
 $L_\o$  and  $d_\o$ denote the 
 Lipschitz constant  and the diamater of $\o$, respectively; moreover, here and in similar occurrences in what follows, the dependence of a constant on $L_\o$  and  $d_\o$ is understood just via an upper bound for them. 
\par A definition of capacity relative to an open set is  recalled in Section \ref{diff}. 
Let  us just mention here that, if $n \geq 3$, then the  capacity ${\rm cap}_{B_1(x)}$ relative to $B_1(x)$ is equivalent to the standard capacity in the whole of $\rn$, up to multiplicative constants depending on $n$.

\begin{theorem}\label{seconddir} {\rm {\bf [Global estimate under capacitary conditions on curvatures]}}
 Let $\Omega$ be a bounded Lipschitz  domain in $\rn$, $n \geq 2$ such that   $\partial \Omega \in  W^{2,1}$.  
  Assume that $p >\tfrac 32$, and  ${\bf f} \in
L^2(\Omega, \rN )$.  Let $\bu$ be an approximable  solution   to either  the
Dirichlet problem \eqref{localeq}$+$\eqref{dircond}, or   the Neumann problem
\eqref{localeq}$+$\eqref{neucond}. There exists a constant $c=c(n, N, p, L_\Omega, d_\o)$ such that, if  $\Omega$ fulfills \eqref{capcond}  for such a constant $c$,
then
\begin{equation}\label{seconddir1} |\nabla \bu|^{p-2} \nabla \bu
\in W^{1,2}(\Omega, \rNn ).
\end{equation}
Moreover,
\begin{equation}\label{seconddir2} C_1 \|{\bf f}\|_{L^2(\Omega, \rN)} \leq \||\nabla \bu|^{p-2} \nabla \bu\|_{W^{1,2}(\Omega, \rNn)} \leq C_2 \|{\bf f}\|_{L^2(\Omega, \rN)}
\end{equation}
for some positive constants $C_1=C_1(n,  N, p)$ and $C_2=C_2(\o, N, p)$.
%
%
\end{theorem}

\begin{remark}\label{sharpcap} {\rm 
Assumption \eqref{capcond} cannot be relaxed by just requring that the limit in \eqref{capcond} be finite. Indeed, examples of domains $\o$ can be produced for which $\mathcal K_\o(r)<\infty$ for $r\in (0,1)$, where the solution to the Dirichlet problem for the Laplace operator, with a suitable smooth right-hand side,   belongs to $W^{2,2}(\o)$ if and only if the limit  in  \eqref{capcond} does not exceed some explicit threshold \cite{Ma67, Ma73}.
%
%
The boundary of the relevant domains is smooth, save a small portion which agrees with the graph of a function $\Theta$ of the variables $(x_1, \dots , x_{n-1})$ having the form $\Theta (x_1, \dots , x_{n-1})= c|x_1|(\log |x_1|)^{-1}$ for some constant $c$ and for small $x_1$.  \color{black}
}
\end{remark}

Capacities can be dispensed in the description of the  boundary regularity for the conclusions of Theorem \ref{seconddir} to hold, and can be replaced with a 
condition involving only integrability properties of $\mathcal B$. The relevant condition reads
\begin{equation}\label{smalln}
 \lim _{r\to 0^+} \Big(\sup _{x \in \partial \Omega} \|\mathcal B \|_{X(\partial \Omega \cap B_r(x))}\Big) < c \,,
\end{equation}
for a suitable constant $c=c(n, N, p, L_\Omega, d_\Omega)$,
where 
\begin{equation}\label{Xspace}
X= \begin{cases} L^{n-1, \infty} & \quad \hbox{if $n \ge 3$,}
\\
L^{1, \infty} \log L  & \quad \hbox{if $n =2$.}
\end{cases}
\end{equation}
The spaces appearing on the right-hand side of \eqref{Xspace} are Marcinkiewicz type spaces, also called weak type space. Precisely,  $L^{n-1, \infty}$ denotes the
weak Lebesgue space $L^{n-1}$, and  $L^{1, \infty} \log L$   the weak Zygmund space $L\log L$. Let us point out that the assumption $\partial \o \in W^2X$ and condition \eqref{smalln} do not even entail that $\partial \o \in C^1$. Here, the notation $\partial \o \in W^2X$ means that the functions which  locally  describe $\partial \o$ are twice weakly differentiable, with second derivatives in  $X$.
\par As observed in Remark \ref{sharpcond} below, condition \eqref{smalln} is minimal as far as integrability properties of the curvatures of $\partial \o$ ar concerned. However, 
 there exist open sets where the capacitary criterion of Theorem \ref{seconddir} applies, whereas the integrability condition \eqref{smalln} assumed in the following corollary does not.

\color{black}

\begin{corollary}\label{secondmarc} {\rm {\bf [Global estimate under integrability conditions on  curvatures]}}
 Let $\Omega$ be a bounded Lipschitz  domain in $\rn$, $n \geq 2$ such that   $\partial \Omega \in  X$. Let $p$, $\bm f$ and $\bm u$ be as in Theorem \ref{seconddir}.
There exists a constant $c=c(n, N, p, L_\Omega, d_\o)$ such that, if  $\Omega$ fulfills \eqref{smalln}  for such a constant $c$,
then $|\nabla \bu|^{p-2} \nabla \bu
\in W^{1,2}(\Omega, \rNn )$, and inequality \eqref{seconddir2} holds
%
%
%
for some positive constants $C_1=C_1(n,  N, p)$ and $C_2=C_2(\o, N, p)$.
%
%
\end{corollary}

\begin{remark}\label{sharpcond}
 {\rm  The sharpness of 
assumption \eqref{smalln} in Corollary  \ref{secondmarc}   can be demonstrated, for instance, when $n=3$ and   $p\in (\tfrac 32,2]$, by an 
 example from \cite{KrolM}.  In that paper,  open sets  $\o \subset \mathbb R^3$,  with $\partial \o \in W^2L^{2, \infty}$ but whose limit in \eqref{smalln} is not small enough, are exhibited where   the  solution $u$ to the Dirichlet problem, with  a smooth right-hand side,  is such that $|\nabla u|^{p-2}\nabla u \notin W^{1,2}(\Omega)$.    If $n=2$ and $p=2$, so that the $p$-Laplace operator reduces to the classical Laplacian, there exist open sets $\Omega$, with 
$\partial \o \in  W^2L^{1, \infty} \log L$, for which the limit in  \eqref{smalln}   is larger than some critical value,
and  where the  solution $u$ to  the Dirichlet problem  with  a smooth right-hand side fails to belong to $W^{2,2}(\o)$ \cite{Ma67} (see also   \cite[Section 14.6.1]{MaSh}, where Neumann problems are discussed). 
\\ Let us notice that the domains mentioned in Remark  \ref{sharpcap} are such that $\partial \o \notin   W^2L^{2, \infty}$ if $n \geq 3$, and hence, unlike Theorem \ref{seconddir}, Corollary  \ref{secondmarc}  cannot be applied. 
}
\end{remark}

\begin{remark}\label{rembound}
 {\rm  
Condition \eqref{smalln}  holds, for instance, if $n \geq 3$ and $\partial \o \in W^{2,n-1}$, and if $n=2$ and $\partial \o \in W^{2}L\log L$ (and hence, a fortiori,  if  $\partial \o \in W^{2,q}$ for some $q>1$).  In  fact, in all these cases  
\begin{equation}\label{limK=0}
\lim _{r \to 0^+}\mathcal K_\o(r)=0\,.
\end{equation}
In particular, condition \eqref{limK=0}, and hence \eqref{smalln}, is satisfied by any domain with a smooth boundary.
%
%
%
}
\end{remark}

\par In the case when $\o$ is   convex, the conclusions of Theorem \ref{seconddir} do not require any additional  assumption on $\partial \o$.   As will be apparent from its proof, such a statement holds thanks to the fact that the curvatures (in any weak sense) on $\partial \o$ have  a  sign that  enables us to disregard certain boundary integrals arising in our estimates. 


\begin{theorem}\label{secondconvex} {\rm {\bf [Global estimate in convex domains]}}
Let $\o$ be any  bounded convex open set in $\rn$, with $n \geq 2$. Let $p$, $\bm f$ and $\bm u$ be as in Theorem \ref{seconddir}.
Then  $|\nabla \bu|^{p-2} \nabla \bu
\in W^{1,2}(\Omega, \rNn )$, and inequality \eqref{seconddir2} holds
for some positive constants $C_1=C_1(n, N, p)$ and $C_2=C_2(\o, N, p)$.
\end{theorem}

A local counterpart of the above theorems   is the subject of  our last main result. In the statement, $B_R$ and $B_{2R}$ denote concentric   balls of radius $R$ and $2R$, respectively.

\begin{theorem}\label{secondloc} {\rm {\bf [Local estimate]}}
Let $\Omega$ be an open set in $\mathbb R^n$, with $n \geq 2$. Assume that $p >\tfrac 32$, and $\bm{f} \in L^2_{\rm loc}(\Omega, \rN)$.
   Let $\bu$ be an approximable  local solution   to  system
\eqref{localeq}.
Then
\begin{equation}\label{secondloc1} |\nabla \bu|^{p-2} \nabla \bu
\in W^{1,2}_{\rm loc}(\Omega, \rNn ),
\end{equation}
and there exists a constant $C=C(n, N, p)$ such that
\begin{equation}\label{secondloc2} \||\nabla \bu|^{p-2} \nabla \bu\|_{W^{1,2}(B_R, \rNn)} \leq C \Big(\|\bm{f}\|_{L^2(B_{2R}, \rN)} + (R^{-\frac n2}+ R^{-\frac n2 -1})\|\nabla \bu\|_{L^{p-1}(B_{2R}, \rNn)}^{p-1} \Big)
\end{equation}
for any ball $B_{2R} \subset \subset \Omega$.
\end{theorem}

\section{Key inequalities}\label{diff}



This section is devoted to the proof of local and global norm inequalities involving 
a smooth differential  operator, modeled upon the $p$-Laplacian. Here, the global results are established in   domains with a smooth boundary. They are the subject of 
Theorem \ref{step1} below, and constitute  a fundamental step in our  approach. 
\\ The differential operators coming into play are built upon  nonnegative functions 
 $a \in C^1([0, \infty))$, satisfying suitable assumptions,  that replace the function $t^{p-2}$ appearing in the $p$-Laplace operator. 
For any such function $a$, we
define
 \begin{equation}\label{ia}
i_a= \inf _{t >0} \frac{t a'(t)}{a(t)} \qquad \hbox{and} \qquad
s_a= \sup _{t >0} \frac{t a'(t)}{a(t)},
\end{equation}
where $a'$ stands for the derivative of $a$.

\begin{theorem}\label{step1}
Let $n \geq 2$, $N \geq 1$, and let $\o$ be an open set in $\rn$. Assume that $a: [0, \infty) \to [0, \infty)$ is a function having the form $a(t)=\widehat a(t^2)$ for some function $\widehat a \in  C^1([0, \infty))$, and such that $a(t)>0$ for $t>0$.  Suppose that
\begin{equation}\label{align}
i_a>-\tfrac 12\,,
\end{equation}
and  
\begin{equation}\label{sup}
s_a< \infty\,.
\end{equation}
(i) There exists a constant $C= C(n,N, i_a, s_a)$,  such that
\begin{multline}\label{loc15}
\|a(|\nabla \bu|)\nabla \bu\|_{W^{1,2}(B_R, \mathbb R^{Nn})} \\ \leq C\big(\|{\rm {\bf div}} ( a(|\nabla \bu|) \nabla \bu)\|_{L^2(B_{2R}, \mathbb R^{N})} + (R^{-\frac n2}+R^{-\frac n2-1})\|a(|\nabla \bu|)\nabla \bu\|_{L^1(B_{2R}, \mathbb R^{Nn})}\big)
\end{multline}
for every function $\bu \in C^3(\o, \rN)$ and any ball $B_{2R}  \subset \subset \o$.
\\ (ii) Assume, in addition, that $\o$ is a bounded open set with $\partial \o \in C^2$.
There exists a constant $c=c(n,N, i_a, s_a, L_\o, d_\o)$ such that, if
\begin{equation}\label{KK}
\mathcal K_\o (r) \leq \mathcal  K(r) \quad \hbox{for $r \in (0,1)$,}
\end{equation}
for some function $\mathcal K : (0,1) \to [0, \infty)$ satisfying
\begin{equation}\label{limK}
\lim _{r \to 0^+} \mathcal  K(r) <c\,,
\end{equation}
then
%
\begin{align}\label{fund}
\|a(|\nabla \bu|) \nabla \bu\|_{W^{1,2}(\o, \rNn)} \leq C \big(\|{\rm {\bf div}} ( a(|\nabla \bu|) \nabla \bu)\|_{L^2(\o, \rN)}+ 
\| a(|\nabla \bu|) \nabla \bu\|_{L^1(\o, \rNn)}\big)
\,
\end{align}
for some constant $C=C(n,N, i_a, s_a, L_\o, d_\o, \mathcal  K)$, and for
 every function $\bu \in C^3(\o, \rN)\cap C^2(\overline \o, \rN)$ fulfilling either \eqref{dircond} or \eqref{neucond}.
%
%
%
\\
In particular, if  $\o$ is convex, then inequality \eqref{fund} holds whatever $\mathcal K_\o$ is, and the constant $C$ in \eqref{fund} only depends on $n,N, i_a, s_a, L_\o, d_\o$. 
\end{theorem}

\begin{remark}\label{constant} {\rm 
By Remark \ref{rembound}, equation \eqref{limK}  holds, with $\mathcal K = \mathcal K _\o$,  for each single open set $\o$ with   $\partial \o \in C^2$. The dependence of  the constant $C$ in inequality \eqref{fund} just through an upper bound $\mathcal K$ for $\mathcal K_\o$, satisfying \eqref{limK},  will be crucial in approximating domains $\o$ with minimal boundary regularity as  in Theorem \ref{seconddir}.}
\end{remark}

The  proof of Theorem \ref{step1} relies upon the differential inequality contained in the next lemma. In what follows, we shall use the notation $\bm u = (u^1, \dots , u^N)$ for a vector-valued function $\bm u : \Omega \to \rN$. Moreover, $\lq\lq \, \cdot \, "$ stands for scalar product, and $|\nabla ^2 \bu| = \big( \sum_{\alpha =1}^N\sum _{i,j=1}^n (u_{x_i x_j}^\alpha)^2)^{\frac 12}$.

%

\begin{lemma}\label{lemma1}
Let $n \geq 2$, $N \geq 1$, and let $\o$ be an open set in $\rn$.  Assume that $a: [0, \infty) \to [0, \infty)$ is  a function having the form $a(t)=\widehat a(t^2)$ for some function $\widehat a \in  C^1([0, \infty))$, and such that $a(t)>0$ for $t>0$.  Suppose that condition
 \eqref{align} is in force.
%
%
Then there exists a positive constant $C= C(n, N, i_a)$ such that 
\begin{align}\label{pointwise}
\big|{\rm {\bf div}} (a(|\nabla \bu|)\nabla \bm u )\big|^2 & \geq \sum _{j=1}^n
\big(a(|\nabla \bu|)^2 \bm u_{x_j} \cdot \Delta \bm u \big)_{x_j} \\ \nonumber &-
\sum _{i=1}^n \Big(a(|\nabla \bu|)^2 \sum _{j=1}^n \bm u_{x_j}  \cdot \bm u_{x_i x_j}\Big)_{x_i}
+ Ca(|\nabla \bu|)^2 |\nabla ^2 \bu|^2 \quad \hbox{in $\o$,}
%
%
%
%
\end{align}
for every function $\bu \in C^3(\Omega, \rN )$.
%
\end{lemma}
\par\noindent
{\bf Proof}.  
Computations show that
\begin{align}\label{point1}
\big|& {\rm {\bf div}}  (a(|\nabla \bu|)\nabla \bu )\big|^2  =
\sum_{\alpha =1}^N \big({\rm  div}  (a(|\nabla \bu|)\nabla u^\alpha )\big)^2
\\ \nonumber
& = \sum_{\alpha =1}^N\big(
a(|\nabla \bu|) \Delta u^\alpha + a'(|\nabla \bu|) \nabla |\nabla \bu| \cdot
\nabla u^\alpha \big)^2 \\ \nonumber & =  a(|\nabla \bu|)^2 \sum_{\alpha =1}^N \big((\Delta u^\alpha)^2
- |\nabla ^2 u^\alpha|^2\big) + a(|\nabla \bu|)^2 \sum_{\alpha =1}^N |\nabla ^2 u^\alpha|^2 +
 \\ \nonumber & \quad + a'(|\nabla \bu|)^2 \sum_{\alpha =1}^N(\nabla |\nabla \bu| \cdot \nabla u^\alpha )^2 + 2
a(|\nabla \bu|) a'(|\nabla \bu|)\sum_{\alpha =1}^N\Delta u^\alpha \nabla |\nabla \bu| \cdot
\nabla u^\alpha \\ \nonumber & =  a(|\nabla \bu|)^2\sum_{\alpha =1}^N\Big(\sum _{j=1}^n
(u_{x_j}^\alpha\Delta u^\alpha)_{x_j} - \sum _{i,j=1}^n (u_{x_j}^\alpha u_{x_ix_j}^\alpha)_{x_i}
\Big) + a(|\nabla \bu|)^2 \sum_{\alpha =1}^N|\nabla ^2 u^\alpha|^2
\\ \nonumber & \quad + a'(|\nabla \bu|)^2 \sum_{\alpha =1}^N(\nabla |\nabla \bu| \cdot \nabla u^\alpha )^2
+ 2 a(|\nabla \bu|) a'(|\nabla \bu|)\sum_{\alpha =1}^N \Delta u^\alpha \nabla |\nabla \bu| \cdot
\nabla u^\alpha
\\ \nonumber & = a(|\nabla \bu|)^2\sum_{\alpha =1}^N\Big(\sum _{j=1}^n
(u_{x_j}^\alpha\Delta u^\alpha)_{x_j} - \sum _{i,j=1}^n (u_{x_j}^\alpha u_{x_ix_j}^\alpha)_{x_i}
\Big) 
\\ \nonumber & \quad  + a(|\nabla \bu|)^2 \sum_{\alpha =1}^N \sum _{i, j=1}^n (u_{x_ix_j}^\alpha)^2 + 
\bigg(\frac{a'(|\nabla \bu|)}{|\nabla \bu|}\bigg)^2 \sum_{\alpha =1}^N \Big( \sum_{\beta =1}^N\sum_{k, j  =1}^n  u_{x_k}^\beta u^{\beta}_{x_kx_j} u_{x_j}^\alpha \Big)^2
\\ \nonumber & \quad  + 2 a(|\nabla \bu|) \frac{a'(|\nabla \bu|)}{|\nabla \bu|}
\sum_{\alpha, \beta =1}^N \sum_{k, j  =1}^n\Delta u^\alpha  u_{x_k}^\beta u^{\beta}_{x_kx_j} u_{x_j}^\alpha\,.
\end{align}
Here,  the expression $\tfrac{a'(|\nabla \bu|)}{|\nabla \bu|}$ has to be interpreted as $2\,\widehat a'(0)$ if $\nabla \bm u=0$. 
One has that 
\begin{align}\label{pointsys1}
\sum_{\alpha =1}^N & \sum _{j=1}^n
(a(|\nabla \bu|)^2 u_{x_j}^\alpha\Delta u^\alpha)_{x_j}
\\ \nonumber &= 
a(|\nabla \bu|)^2  \sum_{\alpha =1}^N \sum _{j=1}^n
(u_{x_j}^\alpha\Delta u^\alpha)_{x_j} + 2 a(|\nabla \bu|) \frac{a'(|\nabla \bu|)}{|\nabla \bu|} \sum_{\alpha ,\beta =1}^N \sum _{j,k=1}^n u_{x_k}^\beta  u^{\beta}_{x_kx_j} u_{x_j}^\alpha \Delta u^\alpha\,,
\end{align}
and 
\begin{align}\label{pointsys2}
 \sum_{\alpha =1}^N& \sum _{i,j=1}^n (a(|\nabla \bu|)^2 u_{x_j}^\alpha
u_{x_ix_j}^\alpha)_{x_i} 
\\ \nonumber & =  a(|\nabla \bu|)^2 \sum_{\alpha =1}^N\sum _{i,j=1}^n (u_{x_j}^\alpha u_{x_ix_j}^\alpha)_{x_i} +  2 a(|\nabla \bu|) \frac{a'(|\nabla \bu|)}{|\nabla \bu|}  \sum_{\alpha, \beta =1}^N\sum _{i,j,k=1}^n  u_{x_k}^\beta  u^{\beta}_{x_kx_i} u_{x_j}^\alpha u^\alpha_{x_ix_j}\,.
\end{align}
On making use of equations \eqref{pointsys1}--\eqref{pointsys2}, the equality between the leftmost  and the rightmost sides of equation \eqref{point1} takes the form
\begin{align}\label{pointsys3neq0}
\big|  {\rm {\bf div}}  (a(|\nabla \bu|)\nabla \bu )\big|^2 &=
   \sum _{j=1}^n
(a(|\nabla \bu|)^2 \bm u_{x_j} \cdot \Delta \bm u )_{x_j}  -    \sum _{i,j=1}^n (a(|\nabla \bu|)^2 \bm u_{x_j} 
\cdot \bm u_{x_ix_j})_{x_i}
\\ \nonumber & \quad + 
\bigg(\frac{a'(|\nabla \bu|) }{|\nabla \bu|}\bigg)^2 
\sum_{\alpha =1}^N \Big( \sum_{\beta =1}^N\sum_{k, j  =1}^n u_{x_k}^\beta    u_{x_j}^\alpha  u^{\beta}_{x_kx_j}\Big)^2
\\ \nonumber &  \quad + 2a(|\nabla \bu|) \frac{a'(|\nabla \bu|)}{ |\nabla \bu|}\sum_{\alpha, \beta =1}^N\sum _{i,j,k=1}^n  u_{x_k}^\beta   u_{x_j}^\alpha  u^\alpha_{x_jx_i} u^{\beta}_{x_kx_i}
+ a(|\nabla \bu|)^2\sum_{\alpha =1}^N \sum _{i, j=1}^n (u_{x_ix_j}^\alpha)^2\,.
\end{align}
If $\nabla \bm u =0$, then equation \eqref{pointsys3neq0} tells us that inequality \eqref{pointwise} holds  with $C=1$. 
\\ Assume next that $\nabla \bm u\neq 0$. 
Then, equation  \eqref{pointsys3neq0}  can be rewritten as
\begin{align}\label{pointsys3}
%
\big|&  {\rm {\bf div}}  (a(|\nabla \bu|)\nabla \bu )\big|^2 =
   \sum _{j=1}^n
(a(|\nabla \bu|)^2 \bm u_{x_j} \cdot \Delta \bm u )_{x_j}  -    \sum _{i,j=1}^n (a(|\nabla \bu|)^2 \bm u_{x_j} 
\cdot \bm u_{x_ix_j})_{x_i}
\\ \nonumber & \quad \qquad \quad \qquad \qquad+ a(|\nabla \bu|)^2 \bigg[ \bigg(\frac{a'(|\nabla \bu|) |\nabla \bu|}{a(|\nabla \bu|)}\bigg)^2 
\sum_{\alpha =1}^N \Big( \sum_{\beta =1}^N\sum_{k, j  =1}^n \frac {u_{x_k}^\beta}{|\nabla \bu|}  \frac{u_{x_j}^\alpha}{|\nabla \bu|}  u^{\beta}_{x_kx_j}\Big)^2
\\ \nonumber & \quad \qquad \quad \qquad + 2 \frac{a'(|\nabla \bu|) |\nabla \bu|}{a(|\nabla \bu|)}\sum_{\alpha, \beta =1}^N\sum _{i,j,k=1}^n \frac {u_{x_k}^\beta}{|\nabla \bu|} \frac{u_{x_j}^\alpha}{|\nabla \bu|} u^\alpha_{x_jx_i} u^{\beta}_{x_kx_i}
+ \sum_{\alpha =1}^N \sum _{i, j=1}^n (u_{x_ix_j}^\alpha)^2\bigg]\,.
\end{align}
Now, define
$$\vartheta  = \frac {|\nabla \bu|a'(|\nabla \bu|)}{a(|\nabla \bu|)},$$
and 
$$\bm \omega ^\alpha  = \frac{\nabla u^\alpha}{|\nabla \bu|},
\quad 
\bm H ^\alpha = \nabla ^2 u^\alpha \qquad \hbox{for $\alpha =1, \dots , N$.} $$
Observe that $\bm \omega ^\alpha \in \mathbb R^n$, and  $\sum _{\alpha =1}^N|\bm \omega ^\alpha |^2 =\sum _{\alpha =1}^N\sum _{i=1}^n (\omega_i^\alpha)^2 =1$, where we have set 
$\bm \omega ^\alpha  = (\omega _1^\alpha, \dots , \omega _n^\alpha)^T$, where $T$ stands for transpose. Moreover, 
$\bm H ^\alpha$ is a symmetric matrix in $\mathbb R^{n^2}$, and, by 
assumption \eqref{align}, 
\begin{align}\label{3/2}
\vartheta  > -\tfrac 12.
\end{align}
With this notations in place, the expression in square brackets on the right-hand side of \eqref{pointsys3} takes the form
\begin{align}\label{point3}
\vartheta ^2 \sum _{\alpha =1}^N   \Big(\sum_{\beta =1}^N \bm H^\beta  \bm \omega ^\beta \cdot \bm \omega ^\alpha\Big)^2 
+ 2\vartheta  
\sum _{\alpha, \beta =1}^N 
\bm H ^\alpha \bm \omega^\alpha  \cdot \bm H^\beta \bm \omega ^\beta  +  \sum _{\alpha =1}^N  {\rm tr }\big((\bm H ^\alpha )^2\big)\,,
\end{align}
where $\lq\lq {\rm tr }"$ denotes  the trace of a matrix. Since the middle term 
in \eqref{point3} equals $\big\|\sum _{\alpha =1}^N 
\bm H ^\alpha \bm\omega ^\alpha\big\|^2$, 
the proof of inequality \eqref{pointwise} is thus reduced to showing that 
%
%
\begin{align}\label{pointsys4}
\vartheta ^2 \sum _{\alpha =1}^N   \Big(\sum_{\beta =1}^N \bm H^\beta  \bm\omega ^\beta \cdot \bm \omega ^\alpha\Big)^2 
+ 2\vartheta
\bigg\|\sum _{\alpha =1}^N 
\bm H ^\alpha \bm\omega ^\alpha\bigg\|^2  +  \sum _{\alpha =1}^N  {\rm tr }\big((\bm H ^\alpha )^2\big) \geq C  \sum _{\alpha =1}^N  {\rm tr }\big((\bm H ^\alpha )^2\big)
\end{align}
for some   positive constant $C=C(n, N,  i_a)$. It $\vartheta \geq 0$, inequality \eqref{pointsys4} holds with $C=1$. Assume next that
\begin{equation}\label{gen5}
-\tfrac 12 < \vartheta < 0\,.
\end{equation}
For each $\alpha = 1, \dots , N$ fix an orthonormal  basis in $\rn$ of eigenvectors $\bm e_1^\alpha, \dots \, \bm e_n^\alpha$ of $\bm H^\alpha $, in which  $\bm H^\alpha $  takes the diagonal form $\bm \Lambda ^\alpha = {\rm
diag}(\lambda _1^\alpha, \dots \lambda _n^\alpha)$, where $\lambda _1^\alpha, \dots \, \lambda _n^\alpha$ stand for the eigenvalues of $\bm H^\alpha $. Let $\overline {\bm \omega} ^\alpha =(\overline\omega _1^\alpha, \dots , \overline \omega _n^\alpha)^T$ be the vector of the components of $\bm \omega^\alpha$ with respect to this basis. Thus, on denoting by $\bm M^\alpha= (m_{ij}^\alpha)$ the orthogonal matrix whose columns are the vectors $\bm e_j$, with $j=1, \dots , n$, one has that
\begin{equation}\label{change}
\bm \omega ^\alpha = \bm M^\alpha \overline {\bm \omega}^\alpha \quad \quad \quad \hbox{for $\alpha =1, \dots , N$,}
\end{equation}
and 
\begin{equation}\label{quadratic}
\bm H^\alpha = \bm M^\alpha \bm \Lambda ^\alpha (\bm M^\alpha )^T \quad \quad \quad \hbox{for $\alpha =1, \dots , N$.}
\end{equation}
Moreover, on setting $\bm D ^\alpha = \bm M^\alpha \bm \Lambda ^\alpha$, with $\bm D ^\alpha = (d_{ij}^\alpha)$, one has that
\begin{equation}\label{product}
d_{ij}^\alpha = m_{ij}^\alpha \lambda _j^\alpha \qquad \hbox{for $i,j=1, \dots n$ and $\alpha =1, \dots \, N$.}
\end{equation}
Clearly, 
\begin{equation}\label{trace}
\sum _{\alpha =1}^N  {\rm tr }\big((\bm H ^\alpha )^2\big) = 
\sum _{\alpha =1}^N  \sum _{i=1}^n (\lambda _i^\alpha)^2 \,.
\end{equation}
On the other hand, by \eqref{change} and \eqref{quadratic},
$$ \bm H^\alpha \bm \omega ^\alpha = \bm M^\alpha \bm \Lambda ^\alpha \overline{\bm \omega}^\alpha \qquad \hbox{for $\alpha =1, \dots \, N$.}
$$
Hence, the following chain holds:
\begin{align}\label{nov1}
\bigg\|\sum _{\alpha =1}^N 
\bm H ^\alpha \bm\omega ^\alpha\bigg\|^2
 & = 
\bigg\|\sum _{\alpha=1}^N 
\bm M^\alpha \bm \Lambda ^\alpha \overline{\bm \omega}^\alpha \bigg\|^2 
 = \sum _{i=1}^n \bigg(\sum _{\alpha =1}^N  \sum _{j=1}^n d_{ij}^\alpha \overline \omega _j^\alpha\bigg) ^2
\\ \nonumber 
&\leq \sum _{i=1}^n \bigg(\sum _{\alpha =1}^N  \sum _{j=1}^n (d_{ij}^\alpha)^2 \bigg) \bigg( \sum _{\alpha =1}^N\sum _{j=1}^n (\overline \omega _j^\alpha)^2\bigg)
= \sum _{\alpha =1}^N   \sum _{i,j=1}^n (d_{ij}^\alpha)^2  
 =  \sum _{\alpha =1}^N  \sum _{i,j=1}^n (m_{ij}^\alpha)^2 (\lambda _j^\alpha)^2
\\ \nonumber 
& =  \sum _{\alpha =1}^N \bigg( \sum _{j=1}^n \Big((\lambda _j^\alpha)^2
 \sum _{i=1}^n(m_{ij}^\alpha)^2\Big)\bigg) =  \sum _{\alpha =1}^N  \sum _{j=1}^n (\lambda _j^\alpha)^2 = \sum _{\alpha =1}^N  {\rm tr }\big((\bm H ^\alpha )^2\big) \,.
\end{align}
Here, we have made use of the fact that, since $\bm M^\alpha$ is an orthogonal matrix for   $\alpha =1, \dots , N$,
$$\sum _{\alpha =1}^N\sum _{j=1}^n  (\overline \omega _j^\alpha)^2 = \sum _{\alpha =1}^N\sum _{j=1}^n (\omega _j^\alpha)^2=1,$$
and 
$$ \sum _{i=1}^n(m_{ij}^\alpha)^2 =1  \qquad \hbox{for   $j=1, \dots , n$ and $\alpha =1, \dots , N$.}
$$
Owing to  \eqref{gen5}, 
inequality  \eqref{nov1} implies that 
\begin{multline}\label{nov2}
\vartheta ^2 \sum _{\alpha =1}^N   \Big(\sum_{\beta =1}^N \bm H^\beta  \bm\omega ^\beta \cdot \bm \omega ^\alpha\Big)^2 
+ 2\vartheta
\bigg\|\sum _{\alpha =1}^N 
\bm H ^\alpha \bm\omega ^\alpha\bigg\|^2 +  \sum _{\alpha =1}^N  {\rm tr }\big((\bm H ^\alpha )^2\big) \\ \geq (2\vartheta +1)\sum _{\alpha =1}^N  {\rm tr }\big((\bm H ^\alpha )^2\big)
\geq (2i_a +1)\sum _{\alpha =1}^N  {\rm tr }\big((\bm H ^\alpha )^2\big),
\end{multline}
whence inequality \eqref{pointsys4} follows with $C= 2 i_a +1$. Altogether, we have shown that inequality \eqref{pointwise} holds with $C=\min \{2i_a +1, 1\}$.
\qed

\smallskip
\par

The next lemma collects Sobolev type inequalities of a  form suitable for our applications in the proof of Theorem \ref{step1}.

\begin{lemma}\label{lemmasobring} Let $\o$ be a bounded Lipschitz domain in $\rn$, $n \geq 2$. Then  for every $\delta >0$, there exists a constant 
$C=C(n, \delta, L_\o, d_\o)$ 
such that
\begin{align}\label{main13}
\int _\Omega v^2\, dx \leq \delta \int _\Omega |\nabla v|^2\, dx + C \bigg(\int _\Omega |v|\, dx \bigg)^2
\end{align}
for every $v \in W^{1,2}(\Omega)$. 
\\ In particular, there exists a constant $C=C(n)$ such that,
if $0<R\leq \sigma < \tau \leq 2R$,  then 
\begin{equation}\label{sobolevring}
\int _{B_\tau \setminus B_\sigma}  v^2\, dx \leq 	C \delta ^2R^2 \int _{B_\tau \setminus B_\sigma}  |\nabla v|^2\, dx + \frac{C}{ \delta^{n}(\tau - \sigma) R^{n-1}} \bigg(\int _{B_\tau \setminus B_\sigma}   |v|\, dx\bigg)^2
\end{equation}
for every $\delta \in (0,1)$,  every $v \in W^{1,2}(B_\tau \setminus B_\sigma)$, and any concentric balls $B_\sigma $ and $B_\tau$.
\end{lemma}

Inequality \eqref{main13} can be established as in \cite[Proof of Theorem 1.4.6/1]{Mabook}. As for inequality \eqref{sobolevring}, see \cite[Proof of Theorem 2.1]{cmsecond}.
\par
A result of \cite{Mazya62,  Mazya64} (see also \cite[Section 2.5.2]{Mabook}) provides us with a necessary and sufficient capacitary condition for the validity of a Poincar\'e type inequality with measure. A special case of that result is enucleated in Lemma \ref{tracecap} below.
Recall that, given an open set $\o \subset \rn$, the capacity ${\rm cap}_\o(E)$  of a set $E  \subset \o$ relative to $\o$ is defined as
\begin{equation}\label{cap}
{\rm cap}_\o(E) = \inf \bigg\{\int _\o |\nabla v|^2\, dx : v\in C^{0,1}_0(\o), v\geq 1 \,\hbox{on }\, E\bigg\}.
\end{equation}
Here, $C^{0,1}_0(\o)$ denotes the space of  Lipschitz continuous, compactly supported functions in $\o$. 
%
\color{black}

\begin{lemma}\label{tracecap}
Assume that $\o$ is the subgraph  in $\rn$,   $n \geq 2$, of a Lipschitz continuous function of $(n-1)$ variables.  Let $\varrho$ be a nonnegative function in  $L^1_{\rm loc}(\partial \o)$, let $x \in \partial \o$, $r_0\in (0,1)$ and $r \in (0, r_0)$.  Then the following properties are equivalent:
\\ (i) There exists a constant $C_1$ such that
\begin{equation}\label{main6cap}
\int _{\partial \Omega \cap B_r(x)} v^2 \, \varrho\, d \hh  \leq C _1
\int _{\Omega \cap B_r(x)}|\nabla v|^2 \, dy
\end{equation}
for every   $v\in C^{0,1}_0(B_r(x))$.
\\ (ii) There exists a constant $C_2$ such that
\begin{equation}\label{main6cap3} \int _E\, \varrho\, d\hh \leq C_2\,{\rm cap}_{B_1(x)}(E),
\end{equation}
for every set $E\subset \partial \o  \cap B_r(x)$.
\\ Moreover, the best 
constans $C_1$ and $C_2$ in \eqref{main6cap} and \eqref{main6cap3}
 are equivalent, up to multiplicative constants depending on $n$, $r_0$ and  on an upper bound for $L_\o$.
\\ The same statement holds if $\o$ is a bounded Lipschitz domain in $\rn$, and  $r_0 \in (0,1)$ is a suitable positive number depending on $n$ and on an upper bound for $d_\o$ and   $L_\o$. In this case, the best constants $C_1$ and $C_2$ are equivalent up to multiplicative constants depending on $n$, $r_0$, $L_\o$ and $d_\o$.
\end{lemma}

\begin{remark}\label{pallepiene}
{\rm A flattening argument for $\partial \o$, combined with an even extension of trial functions, shows that inequality \eqref{main6cap}, and hence also \eqref{main6cap3}, is equivalent to
\begin{equation}\label{main6cap2}
\int _{\partial \Omega \cap B_r(x)} v^2 \, \varrho\, d \hh \leq C_3 
\int _{B_r(x)}|\nabla v|^2 \, dy
\end{equation}
for some constant $C_3$, and every  $v \in C^{0,1}_0(B_r(x))$. Furthermore,  the best constant $C_3$ is equivalent (up to multiplicative constants depending on $n$, $r_0$ and  $L_\o$) to the best constans $C_1$ and $C_2$  in \eqref{main6cap} and   \eqref{main6cap3}.}
\end{remark}

As a consequence of Lemma \ref{tracecap}, a capacity-free condition  for the validity of a weighted trace inequality on balls centered on the boundary can be derived 
 -- see  \cite[Proof of Theorem 2.4]{cmsecond}. The relevant condition can be formulated in terms of membership of the weight in the Marcinkiewicz type spaces on $\partial \o$, with respect to the measure $\hh$, appearing in \eqref{Xspace}. Importantly, the constant in the Poincar\'e inequality turns out to be proportional to the norm of the weight in the Marcinkiewicz   space.
\\ Recall that 
 the Marcinkiewicz space $L^{q, \infty} (\partial \Omega)$, also called weak $L^q(\partial \Omega)$ space, is the Banach function space endowed with the norm defined as 
\begin{equation}\label{weakleb}
\|\psi\|_{L^{q, \infty} (\partial \Omega)} = \sup _{s \in (0, \hh(\partial \Omega))} s ^{\frac 1q} \psi^{**}(s)
\end{equation}
for a measurable function $\psi$ on $\partial \Omega$. Here, $\psi ^{**}(s)= \smallint _0^s \psi^* (r)\, dr$ for $s >0$, where $\psi^*$ denotes the decreasing rearrangement of $\psi$.
%
The Marcinkiewicz type space $L^{1, \infty} \log L (\partial \Omega)$ is equipped with the norm given by
\begin{equation}\label{weaklog}
\|\psi\|_{L^{1, \infty} \log L (\partial \Omega)} = \sup _{s \in (0, \hh(\partial \Omega))} s \log\big(1+ \tfrac{C}s\big) \psi^{**}(s),
\end{equation}
for any constant  $C>\hh(\partial \Omega)$.  Different constants $C$ result in equivalent norms in \eqref{weaklog}.

\begin{lemma}\label{traceineq}
Let $\o$ be a bounded Lipschitz domain in $\rn$, $n \geq 2$.  Let $\varrho$ be nonnegative function in $X(\partial \o)$, where  $X$ is the space appearing in \eqref{Xspace}. Then, there exist  constants  $C>0$ and $r_0\in (0,1)$, depending on $n$ and on an upper bound for $d_\o$ and $L_\o$, such that
\begin{equation}\label{main6X}
\int _{\partial \Omega \cap B_r(x)} v^2 \, \varrho \, d \hh  \leq C \sup_{y \in \partial \o} \|\varrho\|_{X(\partial \o\cap B_r(y))}
\int _{\Omega \cap B_r(x)}|\nabla v|^2 \, dy
\end{equation}
for every $x \in \partial \Omega$, $r\in (0, r_0)$, and 
 every   function $v\in C^{0,1}_0(B_r(x))$.
%
%
%
%
%
\end{lemma}

We are now in a position  to prove Theorem \ref{step1}.

\smallskip
\par\noindent
{\bf Proof of Theorem \ref{step1}}. Assume that $\bu \in C^3(\o, \rN)$, and 
let $\xi \in C^\infty _0(\mathbb R^n)$. Multiplying through   inequality  \eqref{pointwise} by $\xi^2$, and integrating both sides of the resulting inequality over $\Omega$  yield
\begin{align}\label{main1}
\int _\Omega \xi^2  |{\rm  {\bf div}}&( a(|\nabla \bu|) \nabla \bu)|^2\, dx \\ \nonumber &
\geq \int_\Omega \xi ^2 \bigg[
\sum _{j=1}^n
\big(a(|\nabla \bu|)^2 \bm u_{x_j} \cdot \Delta \bm u \big)_{x_j}  -
\sum _{i=1}^n \Big(a(|\nabla \bu|)^2 \sum _{j=1}^n \bm u_{x_j}  \cdot \bm u_{x_i x_j}\Big)_{x_i}\bigg]\, dx
\\ \nonumber & \quad + C\int _\Omega \xi ^2 a(|\nabla \bu|)^2 |\nabla ^2 \bu|^2\, dx
 \\ \nonumber &  = \int_\Omega \xi ^2  \sum_{\alpha =1}^N\bigg[\sum _{j=1}^n
\big(a(|\nabla \bu|)^2 u_{x_j}^\alpha\Delta u^\alpha\big)_{x_j} -
\sum _{i=1}^n \Big(a(|\nabla \bu|)^2 \sum _{j=1}^n u_{x_j}^\alpha u_{x_i x_j}^\alpha \Big)_{x_i}\bigg]\, dx 
\\ \nonumber & \quad + 
C \int _\Omega \xi ^2 a(|\nabla \bu|)^2 |\nabla ^2 \bu|^2\, dx\,
\end{align}
for some constant $C =C (n, N, i_a)$.
\\
Let us first focus on Part (i) of the statement. Assume that $\xi \in C^\infty _0(\o)$. Then the divergence theorem yields 
\begin{align}\label{main2loc}
\int_\Omega \xi ^2&  \sum_{\alpha =1}^N\bigg[\sum _{j=1}^n
\big(a(|\nabla \bu|)^2 u_{x_j}^\alpha\Delta u^\alpha\big)_{x_j} -
\sum _{i=1}^n \Big(a(|\nabla \bu|)^2 \sum _{j=1}^n u_{x_j}^\alpha u_{x_i x_j}^\alpha \Big)_{x_i}\bigg]\, dx 
\\ \nonumber & \quad = - 2 \int _\Omega  a(|\nabla \bu|)^2\xi \nabla \xi  \cdot  \sum_{\alpha =1}^N \bigg[\Delta u^\alpha \nabla u^\alpha - \sum _{j=1}^n u_{x_j}^\alpha \nabla u_{x_j}^\alpha\bigg]\, dx\,.
\end{align}
By Young's inequality, there exists a constant $C'=C'(n, N)$ such that
\begin{align}\label{main4}
2\bigg|\int _\Omega  & a(|\nabla \bu|)^2\xi \nabla \xi  \cdot  \sum_{\alpha =1}^N \bigg[\Delta u^\alpha \nabla u^\alpha - \sum _{j=1}^n u_{x_j}^\alpha \nabla u_{x_j}^\alpha\bigg]\, dx\bigg| 
\\ \nonumber  & \leq \varepsilon C' \int _\Omega \xi^2 a(|\nabla \bu|)^2 |\nabla^2 \bu|^2\,dx + \frac {C'} \varepsilon \int _\Omega |\nabla \xi|^2 a(|\nabla \bu|)^2 |\nabla \bu|^2\,dx\,
\end{align}
for every  $\varepsilon >0$. Equations \eqref{main1}--\eqref{main4} ensure that 
\begin{align}\label{main3quater}
(C-\varepsilon C')\int _\Omega \xi^2 a(|\nabla \bu|)^2 |\nabla^2 \bu|^2\,dx & \leq 
\int _\Omega \xi^2  |{\rm  {\bf div}} ( a(|\nabla \bu|) \nabla \bu)|^2\, dx   \\ \nonumber & \quad + \frac {C'} \varepsilon \int _\Omega |\nabla \xi|^2 a(|\nabla \bu|)^2 |\nabla \bu|^2\,dx \,.
\end{align}
Let $B_{2R}$ be any ball such that $B_{2R} \subset \subset \Omega$, and let $R\leq \sigma < \tau \leq 2R$. An application of inequality \eqref{main3quater}, with $\varepsilon = \tfrac C{2C'}$ and any function $\xi \in C^{\infty}_0(B_\tau)$ such that $0\leq \xi \leq 1$ in $B_\tau $, $\xi =1 $ in $B_\sigma $ and $|\nabla \xi|\leq C/ (\tau -\sigma)$ for some constant $C=C(n)$, tells us that
\begin{align}\label{loc10}
\int _{B_\sigma} a(|\nabla \bu|)^2 |\nabla^2 \bu|^2\,dx   &  \leq 
C \int _{B_{2R}}  |{\rm  {\bf div}} ( a(|\nabla \bu|) \nabla \bu)|^2\, dx  \\ \nonumber & \quad + \frac C{(\tau - \sigma)^2}\int _{B_\tau \setminus B_\sigma} a(|\nabla \bu|)^2 |\nabla \bu|^2\,dx 
\end{align}
for some constant $C=C(n, N, i_a)$.  Inequality \eqref{sobolevring}, applied with
 $\delta = (\tau -\sigma)/R$ and $v= a(|\nabla u|)u_{x_i}^\alpha$, for $i=1 \dots ,n$ and $\alpha = 1, \dots , N$, yields
\begin{multline}\label{sobolevring1}
\frac {1}{(\tau - \sigma)^2}\int _{B_\tau \setminus B_\sigma}  a(|\nabla \bu|)^2 |\nabla \bu|^2\, dx \\ \leq 	C  \int _{B_\tau \setminus B_\sigma}  a(|\nabla \bu|)^2 |\nabla^2 \bu|^2\, dx + \frac{CR}{ (\tau - \sigma) ^{n+3}} \bigg(\int _{B_\tau \setminus B_\sigma}   a(|\nabla \bu|)|\nabla \bu|\, dx\bigg)^2
\end{multline}
for some constant $C=C(n, N, s_a)$. 
Note that here we have made use of  assumption \eqref{sup} to infer that
\begin{equation}\label{feb20'}
|\nabla (a(|\nabla \bu|)u_{x_i}^\alpha)| \leq C \, a(|\nabla \bu|)|\nabla ^2 \bu| \qquad 	\hbox{a.e. in $\o$,}
\end{equation}
for $i=1, \dots , n$, $\alpha = 1, \dots , N$, and for some constant $C=C(n, N, s_a)$.
 From inequalities \eqref{loc10} and \eqref{sobolevring1} one obtains  that
\begin{align}\label{loc11}
\int _{B_\sigma} a(|\nabla \bu|)^2 |\nabla^2 \bu|^2\,dx   \leq 
 C  \int _{B_\tau \setminus B_\sigma}  a(|\nabla \bu|)^2& |\nabla^2 \bu|^2\, dx 
+
C \int _{B_{2R}}  |{\rm  {\bf div}} ( a(|\nabla \bu|) \nabla \bu)|^2\, dx  
\\ \nonumber &   + \frac{CR}{ (\tau - \sigma) ^{n+3}} \bigg(\int _{B_{2R}}   a(|\nabla \bu|)|\nabla \bu|\, dx\bigg)^2
\end{align}
for some constant $C=C(n, N, i_a, s_a)$, whence
\begin{multline}\label{loc12}
\int _{B_\sigma} a(|\nabla \bu|)^2 |\nabla^2 \bu|^2\,dx    \leq 
 \frac {C}{1+C}  \int _{B_\tau }  a(|\nabla \bu|)^2 |\nabla^2 \bu|^2\, dx 
+
C' \int _{B_{2R}}  |{\rm  {\bf div}} ( a(|\nabla \bu|) \nabla \bu)|^2\, dx  
\\ + \frac{C'R}{ (\tau - \sigma) ^{n+3}} \bigg(\int _{B_{2R}} a(|\nabla \bu|)|\nabla \bu|\, dx\bigg)^2
\end{multline}
for positive constants $C=C(n, N, i_a, s_a)$ and $C'=C'(n, N, i_a, s_a)$. 
Via a  customary iteration argument -- see e.g. \cite[Lemma 3.1, Chapter 5]{giaquinta} -- one can deduce  from 
 \eqref{loc12}  that
\begin{align}\label{loc13}
\int _{B_R} a(|\nabla \bu|)^2 |\nabla^2 \bu|^2\,dx   
 & \leq C \int _{B_{2R}}  |{\rm  {\bf div}} ( a(|\nabla \bu|) \nabla \bu)|^2\, dx  
\\ \nonumber & \quad + \frac{C}{R^{n+2}} \bigg(\int _{B_{2R}} a(|\nabla \bu|)|\nabla \bu|\, dx\bigg)^2
\end{align}
for some constant $C=C(n, N, i_a, s_a)$. Moreover, inequality \eqref{main13}, applied with $\Omega = B_1$, $\delta =1$, and   $v= a(|\nabla u|)u_{x_i}^\alpha$, for $i=1 \dots ,n$ and $\alpha = 1, \dots , N$, 
and a scaling argument imply that there exists a constant $C=C(n, N, s_a)$ such that 
 \begin{align}\label{loc14}
\int _{B_R} a(|\nabla \bu|)^2 |\nabla \bu|^2\,dx  &  \leq 
  \int _{B_R}  a(|\nabla \bu|)^2 |\nabla^2 \bu|^2\, dx 
+ \frac{C}{ R^{n}} \bigg(\int _{B_{R}}   a(|\nabla \bu|)|\nabla \bu|\, dx\bigg)^2\,.
\end{align}
Inequality \eqref{loc15} follows from  \eqref{feb20'}, \eqref{loc13} and \eqref{loc14}.
\\
We now consider Part (ii). 
Let $\o$ be a bounded domain with   $\partial \o \in C^2$. Assume that the function  $\bu \in C^3(\o, \rN)\cap C^2(\overline \o, \rN)$ fulfills either \eqref{dircond} or \eqref{neucond}.
Owing to \cite[Equation
(3,1,1,2)]{Grisvard},
\begin{multline}\label{grisv1}
\Delta u^\alpha \frac{\partial u^\alpha}{\partial \bm \nu } - \sum _{i,j=1}^n
u_{x_i x_j}^\alpha u_{x_i}^\alpha  \nu _j \\ = {\rm div }_T
\bigg(\frac{\partial u^\alpha}{\partial \bm \nu } \nabla _T u^\alpha\bigg) -
{\rm tr}\mathcal B \bigg(\frac{\partial u^\alpha}{\partial \bm \nu }\bigg)^2
- \mathcal B (\nabla _T \, u^\alpha, \nabla _T \, u^\alpha) - 2 \nabla _T\,
u^\alpha \cdot \nabla _T\, \frac{\partial u^\alpha}{\partial \bm \nu }
 \qquad
\hbox{on $\partial \Omega$,}
\end{multline}
 for $\alpha = 1, \dots , N$,
where  ${\rm tr}\mathcal B $ is the trace of $\mathcal B$,  
 ${\rm div }_T$ and
$\nabla _T $ denote the divergence and the gradient
operator on $\partial \Omega$, respectively, and $\nu _j$ stands for the $j$-th component of $\bm \nu$. From 
 the divergence theorem and equation \eqref{grisv1} we deduce that 
\begin{align}\label{main2}
\int_\Omega \xi ^2&  \sum_{\alpha =1}^N\bigg[\sum _{j=1}^n
\big(a(|\nabla \bu|)^2 u_{x_j}^\alpha\Delta u^\alpha\big)_{x_j} -
\sum _{i=1}^n \Big(a(|\nabla \bu|)^2 \sum _{j=1}^n u_{x_j}^\alpha u_{x_i x_j}^\alpha \Big)_{x_i}\bigg]\, dx 
\\ \nonumber & = \int _{\partial \Omega } \xi^2 a(|\nabla \bu|)^2  \sum_{\alpha =1}^N\bigg[\Delta u^\alpha \frac {\partial u^\alpha}{\partial \bm \nu} - \sum _{i,j=1}^nu_{x_i x_j}^\alpha u_{x_i}^\alpha\nu _j \bigg]\, d\mathcal H^{n-1}
\\ \nonumber & \quad  \quad- 2 \int _\Omega  a(|\nabla \bu|)^2\xi \nabla \xi  \cdot  \sum_{\alpha =1}^N \bigg[\Delta u^\alpha \nabla u^\alpha - \sum _{j=1}^n u_{x_j}^\alpha\nabla u_{x_j}^\alpha \bigg]\, dx
\\ \nonumber & = \int _{\partial \Omega } \xi^2 a(|\nabla \bu|)^2  \sum_{\alpha =1}^N \bigg[{\rm div }_T
\bigg(\frac{\partial u^\alpha}{\partial \bm \nu } \nabla _T u^\alpha\bigg) -
{\rm tr}\mathcal B \bigg(\frac{\partial u^\alpha}{\partial \bm \nu }\bigg)^2
\\ \nonumber & \quad \quad  \qquad \qquad \qquad 
- \mathcal B (\nabla _T \, u^\alpha , \nabla _T \, u^\alpha) - 2 \nabla _T\,
u ^\alpha\cdot \nabla _T\, \frac{\partial u^\alpha}{\partial \bm \nu }\bigg] \,d\mathcal H^{n-1}
\\ \nonumber & \quad - 2 \int _\Omega  a(|\nabla \bu|)^2\xi \nabla \xi  \cdot  \sum_{\alpha =1}^N \bigg[\Delta u^\alpha \nabla u^\alpha - \sum _{j=1}^n u_{x_j}^\alpha \nabla u_{x_j}^\alpha\bigg]\, dx\,.
\end{align}
Equations \eqref{main1}, \eqref{main4} and \eqref{main2}  ensure that there exist constants $C=C(n,N, i_a)$ and $C'=C'(n, N)$ such that
\begin{align}\label{main3bis}
(C-\varepsilon C') \int _\Omega \xi^2 a(|\nabla \bu|)^2 |\nabla^2 \bu|^2\,dx & \leq 
\int _\Omega \xi^2  |{\rm  {\bf div}} ( a(|\nabla \bu|) \nabla \bu)|^2\, dx   + \frac {C'} \varepsilon \int _\Omega |\nabla \xi|^2 a(|\nabla \bu|)^2 |\nabla \bu|^2\,dx 
\\ \nonumber & \quad + 
\bigg|\int _{\partial \Omega } \xi^2 a(|\nabla \bu|)^2  \sum_{\alpha =1}^N \bigg[{\rm div }_T
\bigg(\frac{\partial u^\alpha}{\partial \bm \nu } \nabla _T u^\alpha\bigg) -
{\rm tr}\mathcal B \bigg(\frac{\partial u^\alpha}{\partial \bm \nu }\bigg)^2
\\ \nonumber & \quad \quad \quad
- \mathcal B (\nabla _T \, u^\alpha , \nabla _T \, u^\alpha) - 2 \nabla _T\,
u^\alpha \cdot \nabla _T\, \frac{\partial u^\alpha}{\partial \bm \nu }\bigg] \,d\mathcal H^{n-1}\bigg|\,.
\end{align}
If $\bu$ satisfies the boundary condition \eqref{dircond}, then $\nabla _T\, u^\alpha=0$ on $\partial \Omega$ for $\alpha =1, \dots , N$, and hence  
\begin{align}\label{main3dir}
\int _{\partial \Omega }& \xi^2 a(|\nabla \bu|)^2  \sum_{\alpha =1}^N  \bigg[{\rm div }_T
\bigg(\frac{\partial u^\alpha}{\partial \bm \nu } \nabla _T \,u^\alpha\bigg) -
{\rm tr}\mathcal B \bigg(\frac{\partial u^\alpha}{\partial \bm \nu }\bigg)^2
\\ \nonumber & \quad \quad \quad \quad \quad \quad
- \mathcal B (\nabla _T \, u^\alpha , \nabla _T \, u^\alpha) - 2 \nabla _T\,
u^\alpha \cdot \nabla _T\, \frac{\partial u^\alpha}{\partial \bm \nu }\bigg] \,d\mathcal H^{n-1}
\\ \nonumber & = - 
\int _{\partial \Omega } \xi^2 a(|\nabla \bu|)^2 
{\rm tr}\mathcal B \sum_{\alpha =1}^N   \bigg(\frac{\partial u^\alpha}{\partial \bm \nu }\bigg)^2 \,d\mathcal H^{n-1}\,.
\end{align}
On the orther hand, if the boundary condition \eqref{neucond} is in force, then $\displaystyle \frac{\partial u^\alpha}{\partial \bm \nu }=0$  on $\partial \Omega$ for $\alpha =1, \dots , N$. Therefore,
\begin{align}\label{main3neu}
\int _{\partial \Omega }& \xi^2 a(|\nabla \bu|)^2  \sum_{\alpha =1}^N  \bigg[{\rm div }_T
\bigg(\frac{\partial u^\alpha}{\partial \bm \nu } \nabla _T \,u^\alpha\bigg) -
{\rm tr}\mathcal B \bigg(\frac{\partial u^\alpha}{\partial \bm \nu }\bigg)^2
\\ \nonumber & \quad \quad \quad \quad \quad \quad
- \mathcal B (\nabla _T \, u^\alpha , \nabla _T \, u^\alpha) - 2 \nabla _T\,
u^\alpha \cdot \nabla _T\, \frac{\partial u^\alpha}{\partial \bm \nu }\bigg] \,d\mathcal H^{n-1}
\\ \nonumber & = - \int _{\partial \Omega } \xi^2 a(|\nabla \bu|)^2  \sum _{\alpha =1}^N\mathcal B (\nabla _T \, u^\alpha , \nabla _T \, u^\alpha)
 \,d\mathcal H^{n-1}\,.
\end{align}
In both cases, one has that 
\begin{align}\label{main3}
\bigg|\int _{\partial \Omega }& \xi^2 a(|\nabla \bu|)^2  \sum_{\alpha =1}^N  \bigg[{\rm div }_T
\bigg(\frac{\partial u^\alpha}{\partial \bm \nu } \nabla _T \,u^\alpha\bigg) -
{\rm tr}\mathcal B \bigg(\frac{\partial u^\alpha}{\partial \bm \nu }\bigg)^2
\\ \nonumber & \quad \quad \quad \quad \quad \quad
- \mathcal B (\nabla _T \, u^\alpha , \nabla _T \, u^\alpha) - 2 \nabla _T\,
u^\alpha \cdot \nabla _T\, \frac{\partial u^\alpha}{\partial \bm \nu }\bigg] \,d\mathcal H^{n-1} \bigg|
\\ \nonumber & \leq C \int _{\partial \Omega }\xi ^2a(|\nabla \bu|)^2 |\nabla \bu|^2 |\mathcal B| \,d\mathcal H^{n-1}\,,
\end{align}
for some constant $C(n,N)$. 
Now, assume that
\begin{equation}\label{feb22}
\xi \in C^\infty _0(B_r(x_0))
\end{equation}
for some $x_0 \in \overline \Omega$ and $r\in (0, r_0)$, where $r_0\in (0,1)$ is as in the last assertion of the statement of Lemma \ref{tracecap}.\\
First, suppose that
$x_0 \in \partial \o$. 
Then an application of inequality \eqref{main6cap}, with $\varrho = |\mathcal B|$ and $v = \xi \,a(|\nabla \bu|) u_{x_i}^\alpha$ , for $i=1, \dots n$ and $\alpha = 1, \dots , N$, yields
\begin{multline}\label{main7}
 \int _{\partial \Omega } \xi^2 \, a(|\nabla \bu|)^2 |\nabla \bu|^2 |\mathcal B| \,d\mathcal H^{n-1}  \\ \leq C \mathcal K_\o(r)\bigg(\int _\Omega  \xi ^2 a(|\nabla \bu|)^2 |\nabla ^2 \bu|^2\, dx  + \int _\Omega  |\nabla \xi |^2 a(|\nabla \bu|)^2 |\nabla \bu|^2\, dx \bigg)
\end{multline}
for some constant $C=C(n, N, s_a, L_\o, d_\o)$. Note that inequality \eqref{feb20'} has also been exploited in deriving \eqref{main7}.
Combining equations  \eqref{main3bis} and \eqref{main7} tells us that 
\begin{multline}\label{main10}
\Big[C_1- \varepsilon C_3 - C_2 \, \mathcal K_\o(r)\Big] \int _\Omega \xi ^2 a(|\nabla \bu|)^2 |\nabla ^2 \bu|^2\, dx \\ \leq \int _\Omega \xi^2  |{\rm  {\bf div}} ( a(|\nabla \bu|) \nabla \bu)|^2\, dx 
+ \bigg[C_2 \, \mathcal K_\o(r) + \frac {C_3}\varepsilon \bigg] \int _\Omega  |\nabla \xi |^2 a(|\nabla \bu|)^2 |\nabla \bu|^2\, dx
\end{multline}
for some constants $C_1=C_1(n, N,  i_a)$, $C_2=C_2(n, N, s_a, L_\o, d_\o)$ and $C_3=C_3(n, N)$. 
%
%
%
%
If condition \eqref{limK} is fulfilled with $c= \tfrac{C_1}{C_2}$, then there exist 
$\ep>0$ and $C>0$, depending on $\o$ only through $L_\o$, $d_\o$, and   $r'\in (0, r_0)$ depending also on the function $\mathcal K$, such that
\begin{equation}\label{gen2}
C_1- \varepsilon C_3 - C_2 \, \mathcal K_\o(r) \geq C_1- \varepsilon C_3 - C_2 \, \mathcal K (r) \geq C
\end{equation}
%
%
provided that  $r \in (0,r']$. 
Therefore, by inequality \eqref{main10},
\begin{align}\label{main11}
\int _\Omega \xi ^2 a(|\nabla \bu|)^2 |\nabla ^2 \bu|^2\, dx & \leq C \int _\Omega \xi^2  |{\rm  {\bf div}} ( a(|\nabla \bu|) \nabla \bu)|^2\, dx 
+  C \int _\Omega  |\nabla \xi |^2 a(|\nabla \bu|)^2 |\nabla \bu|^2\, dx
\end{align}
for some constant $C=C(n, N, i_a, s_a, L_\o, d_\o)$, if  $r \in (0,r']$ in \eqref{feb22}.
\\ When $x_0\in \Omega$ and $B_r(x_0) \subset \o$, inequality \eqref{main11}   just follows from \eqref{main3bis}, inasmuch as the integral over $\partial \Omega$ in the latter inequality vanishes. Also, the constant $C$ in \eqref{main11}  is independent of $L_\Omega$ and $d_\Omega$.
%
%
\\ Now,  there exists $r^{''}\in (0, r')$,   hence depending on $L_\o$, $d_\o$, and   on the function $\mathcal K$, such that   $\overline \o$ admits a finite covering 
 $\{B_{r_k}\}$   by balls $B_{r_k}$, with  $r^{''} \leq r_k \leq r'$, and having the property that either  $B_{r_k}$ is centered on $\partial \o$, or $B_{r_k}\subset   \o$. Note that this covering can be chosen in such a way that the multiplicity of overlapping of the balls $B_{r_k}$ only depends on $n$. Let 
$\{\xi _k\}$ be a family of functions such that $\xi_k \in C^\infty _0(B_{r_k})$ and 
$\{\xi_k^2\}$ is a partition of unity of $\overline \o$ associated with the covering $\{B_{r_k}\}$. 
Thus  $\sum _{k} \xi _k^2 = 1$ in $\overline \o$.  The functions $\xi_k$ can also be chosen so that $|\nabla \xi_k| \leq C/r_k \leq C/r^{''}$ for some absolute constant $C$. Making use of  inequality \eqref{main11} with $\xi = \xi _k$ for each $k$, and adding the resulting inequalities one obtains that
\begin{align}\label{main12}
\int _\Omega  a(|\nabla \bu|)^2 |\nabla ^2 \bu|^2\, dx & \leq C \int _\Omega   |{\rm  {\bf div}} ( a(|\nabla \bu|) \nabla \bu)|^2
\, dx 
+  C \int _\Omega  a(|\nabla \bu|)^2 |\nabla \bu|^2\, dx
\end{align}
for some constant $C= C(n, N,  i_a, s_a, L_\o, d_\o, \mathcal K)$. 
 Given $\delta >0$, an application of inequality \eqref{main13} with $v = a(|\nabla \bu|) u_{x_i}^\alpha$, $i=1, \dots , n$, $\alpha =1, \dots , N$, and equation \eqref{feb20'} tell  us that
\begin{align}\label{main14}
\int _\Omega  a(|\nabla \bu|)^2 |\nabla \bu|^2\, dx \leq \delta C_1 \int _\Omega a(|\nabla \bu|)^2 |\nabla ^2 \bu|^2\, dx + C_2 \bigg(\int _\Omega a(|\nabla \bu|)|\nabla \bu|\, dx \bigg)^2
\end{align}
for some constants $C_1=C_1(n, N, s_a)$ and $C_2= C_2(n,  N, s_a, \delta, L_\o, d_\o )$.
%
%
%
Choose $\sigma = \tfrac{1}{2CC_1}$ in inequality \eqref{main14}, where $C$ is the constant appearing in \eqref{main12}. Coupling the resulting inequality with inequality   \eqref{main12} enables us to conclude that
\begin{align}\label{main16}
 \int _\Omega  a(|\nabla \bu|)^2 |\nabla ^2\bu|^2\, dx 
\leq C
\int _\Omega |{\rm  {\bf div}} ( a(|\nabla \bu|) \nabla \bu)|^2\, dx  +  C
 \bigg(\int _\Omega a(|\nabla \bu|)|\nabla \bu|\, dx \bigg)^2
\end{align}
for some constant $C=C(n,  N, i_a , s_a, L_\o, d_\o, \mathcal K)$. 
Inequalities \eqref{main14} and \eqref{main16} imply  \eqref{fund}, via \eqref{feb20'}. 
\par Finally, under the assumption  that $\o$ is convex, 
the right-hand side of identity \eqref{grisv1} is nonnegative, since, in this case,
$$     - {\rm tr}\mathcal B \sum_{\alpha =1}^N\bigg(\frac {\partial u^\alpha}{\partial \bm \nu}\bigg)^2  \geq 0 \quad  \hbox{and}  \quad 
- \sum_{\alpha =1}^N \mathcal B (\nabla _T \,u^\alpha , \nabla _T \,u^\alpha) \geq 0 \qquad \quad
\hbox{on $\partial \Omega$ \,.}$$
Therefore, inequality \eqref{main3bis} can be replaced with the stronger inequality
\begin{align}\label{main3ter}
(C-\varepsilon C') \int _\Omega \xi^2 a(|\nabla \bu|)^2 |\nabla^2 \bu|^2\,dx & \leq 
\int _\Omega \xi^2  |{\rm  {\bf div}} ( a(|\nabla \bu|) \nabla \bu)|^2\, dx   + \frac {C'} \varepsilon \int _\Omega |\nabla \xi|^2 a(|\nabla \bu|)^2 |\nabla \bu|^2\,dx 
\end{align}
for some constants $C=C(n,N,i_a)$ and $C'=C'(n,N,i_a)$.
 Starting from this inequality, instead of \eqref{main3bis},
estimate \eqref{fund} follows analogously (and even more easily). In particular, the function $\mathcal K_\o$ does not come into play in this case.
\qed

\section{Local estimates}\label{proofsloc}

The proof of Theorem \ref{secondloc} is accomplished in this section. Approximable local solutions to equation \eqref{localeq} considered in its statement can be defined as follows.
\par 
Assume that $\bm f \in L^1_{\rm loc}(\Omega, \rN) \cap (W^{1,p}_0(\Omega, \rN))'$. A function $\bu \in W^{1,p}_{\rm loc}(\o, \mathbb R^N)$ is called a local weak solution to system \eqref{localeq} if
\begin{equation}\label{weaklocl}
\int _{\o'} |\nabla \bu|^{p-2} \nabla \bu \cdot \nabla  \bm{\varphi} \color{black} \,dx = \int _{\Omega '} {\bf f}\cdot \bm{\varphi}  \,dx
\end{equation}
for every open set $\Omega ' \subset \subset \o$, and every function $\bm{\varphi}  \in W^{1,p}_0(\Omega', \mathbb R^N)$.
\\ Assume that $\bm f \in L^q_{\rm loc}(\o, \rN)$ for some $q\geq 1$. 
A function $\bu \in W ^{1,1}_{\rm loc}(\o, \rN)$ is called a local  approximable  solution to system \eqref{localeq} if $\nabla \bu \in L^{p-1}_{\rm loc}(\o)$, 
and
there exists a sequence $\{\bm{f}_k\}\subset  C^\infty (\o, \rN)$, with
$\bm{f}_k \to \bm{f}$   in $L^q_{\rm loc} (\o, \rN)$,
such that the corresponding sequence of local weak solutions $\{\bu_k\}$ to the system
\begin{equation}\label{localeqk}
- {\rm {\bf  div}} (|\nabla \bu_k|^{p-2}\nabla {\bf u} _k) = {\bf f}_k  \quad {\rm in}\,\,\, \o\,,
\end{equation}
satisfies
\begin{equation}\label{approxuloc}
\bu_k \to \bu \quad \hbox{and} \quad \nabla \bu_k \to \nabla \bu \quad \hbox{a.e. in $\o$,}
\end{equation}
and
\begin{equation}\label{approxaloc}
\lim _{k \to \infty} \int_{\o'} |\nabla \bu_k|^{p-1} \, dx =\int_{\o'} |\nabla \bu|^{p-1}\, dx\,
\end{equation}
for every open set $\o' \subset \subset \o$.

\smallskip
\par\noindent
{\bf Proof of Theorem \ref{secondloc}}.  
Assume, for the time being, that 
\begin{equation}\label{fsmoothloc}
\bm{f} \in C^\infty (\Omega, \rN)\,.
\end{equation}
Given $\varepsilon \in (0,1)$, define the function $b_\varepsilon  : \mathbb R \to (0,\infty)$ as
\begin{equation}\label{aeps}
b_\varepsilon (t) = \sqrt {t^2 + \varepsilon^2}\quad \hbox{for $t \in \mathbb R$.}
\end{equation}
Observe that
\begin{equation}\label{epsindex}
i_{b_\ep^{p-2}} = \min \{p-2, 0\} \quad \hbox{and} \quad s_{b_\ep^{p-2}} = \max \{p-2, 0\}\,,
\end{equation}
where $i_{b_\ep^{p-2}}$ and $s_{b_\ep^{p-2}}$ are defined as in \eqref{ia}, with $a$ replaced by $b_\ep^{p-2}$.
Hence,
\begin{equation}\label{epsindex1}
i_{b_\ep^{p-2}} > - \tfrac 12\,,
\end{equation}
inasuch as we are assuming that  $p> \tfrac 32$.
\\
Let $B_{2R} \subset \subset \Omega$. Given  a local weak solution  $\bu$ to system \eqref{localeq},   let  $\bu_\varepsilon \in \bm u + W^{1,p}_0(B_{2R}, \rN)$ be the weak solution to the Dirichlet problem 
\begin{equation}\label{eqeploc}
\begin{cases}
- {\rm{\bf div}} (b_\ep(|\nabla \bu_\ep|)^{p-2}\nabla \bu_\ep ) = \bm{f}  & {\rm in}\,\,\, B_{2R} \\
 \bu_\ep =\bu  &
{\rm on}\,\,\,
\partial B_{2R} \,.
\end{cases}
\end{equation}
To begin with, we notice that 
\begin{align}\label{cinf}
\bu _\ep \in C^{\infty}(B_{2R}, \rN).
\end{align} 
Indeed, by \cite[Corollary 1.26]{sebastian},    $\nabla \bu_\ep \in L^{\infty}_{\rm loc}(B_{2R}, \rNn)$. Hence, via 
\cite[Corollary 1.26 and equation (1.8d) ]{sebastian}, there exists $\alpha \in (0,1)$ such that $\nabla \bu_\ep \in C^{\alpha}_{\rm loc}(B_{2R}, \rNn)$. In particular, 
   $b_\ep (|\nabla \bu _\ep|)^{p-2} \in C^{1,\alpha}_{\rm loc}(B_{2R})$. The  Schauder theory for linear elliptic systems then entails that $\bu _\ep \in C^{2,\alpha}_{\rm loc}(B_{2R}, \rN)$.  Property \eqref{cinf} then follows by iteration, via the Schauder theory again.
\\
Next, let us observe that
there exists a constant $C=C(n, N, p, R)$  such that
\begin{align}\label{cialpha1}
\int _{B_{2R}} |\nabla \bu _\ep|^p\, dx \leq C \bigg(\int _{B_{2R}}|\bm{f}|^{p'}\,dx +  \int _{B_{2R}}|\nabla \bu |^p\, dx + \ep^p\bigg)\,
\end{align}
for $\ep \in (0,1)$. Here, $p'=\tfrac{p}{p-1}$, the H\"older conjugate.
Actually, the use of  $\bu _\ep - \bu \in W^{1,p}_0(B_{2R}, \rn)$ as a test function in the weak formulation of problem \eqref{eqeploc}  yields
\begin{align}\label{cialpha2}
\int _{B_{2R}} b_\ep (|\nabla \bu _\ep|)^{p-2}\nabla \bu _\ep \cdot (\nabla \bu _\ep - \nabla \bu)\, dx = \int _{B_{2R}} \bm{f} \cdot (\bu _\ep - \bu)\, dx\,.
\end{align}
One can verify that there exist  nonnegative constants $c_1=c_1(p)$ and $c_2=c_2(p)$ such that
\begin{equation}\label{equiv}
t^p -  c_1 \ep ^p\leq b_\ep (t)^{p-2} t^2 \leq c_2 (t^p + \ep^p) \quad \hbox{for $t \geq 0$.}
\end{equation}
Thus, equation \eqref{cialpha2}, Young's inequality and Poincar\'e's inequality imply that, given $\delta \in (0,1)$,
\begin{align}\label{cialpha3}
\int _{B_{2R}}|\nabla \bu _\ep|^p\, dx & \leq  \int _{B_{2R}} |\bm{f}|  |\bu _\ep - \bu|\, dx + C \int _{B_{2R}}(|\nabla \bu _\ep|^{p-1} + \ep^{p-1})|\nabla \bu|\, dx   + C R^n \ep ^p
\\ \nonumber & 
\leq    C' \int _{B_{2R}} |\bm{f}|^{p'}\, dx  +  \delta \int _{B_{2R}}   |\bu _\ep - \bu|^p\, dx 
\\ \nonumber & \quad + \delta \int _{B_{2R}}|\nabla \bu _\ep|^p\, dx + C'  \int _{B_{2R}}|\nabla \bu|^p\, dx + C''  \ep ^p
\\ \nonumber & 
\leq    C' \int _{B_{2R}} |\bm{f}|^{p'}\, dx  +    \delta C''' \int _{B_{2R}}|\nabla \bu _\ep|^p\, dx + C'  \int _{B_{2R}}|\nabla \bu|^p\, dx + C''  \ep ^p
\end{align}
for suitable constants $C=C(p)$, $C'=C'(p,\delta)$, $C''=C''(n,N, p, R)$ and $C'''=C'''(n,p,N,R)$, and for $\ep \in (0,1)$. Choosing $\delta$ sufficiently small yields \eqref{cialpha1}. 
\\ Now,
we claim that 
\begin{equation}\label{akm4}
\nabla \bu_\ep \to  \nabla \bu \quad  \hbox{in $L^p(B_{2R}, \rNn)$}
\end{equation}
as $\ep \to 0^+$.
This claim can be verified via the following argument from \cite[Proof of Theorem 4.1]{AKM}. 
Making use of $\bu_\ep - \bu$ as a test function in the weak formulation of systems \eqref{localeq} and \eqref{eqeploc}, subtracting the resulting equalities, and applying H\"older's inequality yield
\begin{align}\label{akm1}
\int_{B_{2R}}&\big(b_\ep(|\nabla \bu _\ep|)^{p-2}\nabla \bu_\ep - b_\ep(|\nabla \bu|)^{p-2}\nabla \bu\big)\cdot (\nabla \bu _\ep - \nabla \bu)\, dx \\ \nonumber & = 
\int_{B_{2R}}\big(|\nabla \bu|^{p-2}\nabla \bu - b_\ep(|\nabla \bu|)^{p-2}\nabla \bu\big) \cdot(\nabla \bu _\ep - \nabla \bu)\, dx 
\\ \nonumber & \leq C\bigg(\int_{B_{2R}} |\nabla \bu_\ep|^p + |\nabla \bu|^p \, dx \bigg)^{\frac 1p}
\bigg(\int_{B_{2R}}  \big||\nabla \bu|^{p-2}\nabla \bu - b_\ep(|\nabla \bu|)^{p-2}\nabla \bu\big|^{p'}\bigg)^{\frac 1{p'}}\,
\end{align}
for some constant $C=C(p)$.
On the other hand, 
\begin{align}\label{V}
C(|\bm{\xi}|^2 + |\bm{\eta}|^2 + \ep^2)^{\frac {p-2}2} |\bm{\xi} - \bm{\eta}|^2  & \leq \big|b_\ep(|\bm{\xi}|)^{\frac {p-2}2} \bm{\xi}  - b_\ep(|\bm{\eta}|)^{\frac {p-2}2} \bm{\eta}\big|^2\\ \nonumber & \leq C'  \big(b_\ep(|\bm{\xi}|)^{p-2}\bm{\xi} - b_\ep(|\bm{\eta}|)^{p-2}\bm{\eta}\big)\cdot (\bm{\xi} - \bm{\eta}) \qquad \hbox{for $\bm{\xi} , \bm{\eta} \in \rNn$,}
\end{align}
 for suitable positive constants $C=C(p)$ and $C'=C'(p)$.
From \eqref{akm1} and \eqref{V} we deduce that there exists a constant $C=C(p)$ such that
\begin{align}\label{akm2}
\int_{B_{2R}}&(|\nabla \bu_\ep|^2 + |\nabla \bu|^2 + \ep^2)^{\frac {p-2}2} |\nabla \bu_\ep - \nabla \bu|^2   \,dx
\\ \nonumber & \leq C\bigg(\int_{B_{2R}} |\nabla \bu_\ep|^p + |\nabla \bu|^p \, dx \bigg)^{\frac 1p}
\bigg(\int_{B_{2R}}  \big||\nabla \bu|^{p-2}\nabla \bu - b_\ep(|\nabla \bu|)^{p-2}\nabla \bu\big|^{p'}\,dx\bigg)^{\frac 1{p'}}\,.
\end{align}
By inequality \eqref{cialpha1}, the first integral on the right-hand side of inequality \eqref{akm2} is uniformly bounded for $\ep \in (0,1)$. Moreover, the second integral converges to $0$ as $\ep \to 0$, by dominated convergence. Thus, 
\begin{equation}\label{akm3}
\lim _{\ep \to 0^+} \int_{B_{2R}}(|\nabla \bu_\ep|^2 + |\nabla \bu|^2 + \ep^2)^{\frac {p-2}2} |\nabla \bu_\ep - \nabla \bu|^2   \,dx =0\,.
%
\end{equation}
If $p \geq 2$, equation \eqref{akm4} follows  from \eqref{akm3}. 
If $1<p<2$, conclusion \eqref{akm4} still holds, since, by H\"older's inequality,
\begin{multline}\label{akm5}
\int_{B_{2R}}|\nabla \bu_\ep - \nabla \bu|^p   \,dx \\ \leq \bigg(\int_{B_{2R}}(|\nabla \bu_\ep|^2 + |\nabla \bu|^2 + \ep^2)^{\frac {p-2}2} |\nabla \bu_\ep - \nabla \bu|^2   \,dx\bigg)^{\frac p2} 
 \bigg(\int_{B_{2R}}(|\nabla \bu_\ep|^2 + |\nabla \bu|^2 + \ep^2)^{\frac {p}2}    \,dx\bigg)^{\frac {2-p}2} \,.
\end{multline}
From property \eqref{akm4} one can deduce that
\begin{equation}\label{loc17}
\|b_\ep(|\nabla \bu_\ep|)^{p-2}\nabla \bu_\ep\|_{L^1(B_{2R}, \rNn)}\leq C
\end{equation}
for some constant $C$ independent of $\varepsilon \in (0,1)$. Owing to \eqref{epsindex1}, the function $b_\ep^{p-2}$ satisfies  the hypotheses   on the function $a$ in Theorem \ref{step1}. An application of inequality \eqref{loc15} of this theorem to the function $\bm u_\ep$,    and the equation in \eqref{eqeploc}, tell us that
\begin{align}\label{loc18}
\|b_\ep(|\nabla \bu_\ep|)^{p-2}\nabla \bu_\ep & \|_{W^{1,2}(B_R, \rNn)}  \\ \nonumber  & \leq C\big(\|\bm{f}\|_{L^2(B_{2R}, \rN)} + (R^{-\frac n2}+R^{-\frac n2-1})\|b_\ep(|\nabla \bu_\ep|)^{p-2}\nabla \bu_\ep\|_{L^1(B_{2R}, \rNn)}\big)\,,
\end{align}
where     $C=C(n,N, p)$, and, in particular, is indepedent of $\ep$. Inequalities \eqref{loc17} and \eqref{loc18} ensure that the sequence $\{b_\ep(|\nabla \bu_\ep|)^{p-2}\nabla \bu_\ep\}$ is bounded in $W^{1,2}(B_R, \rNn)$, and hence there exists a function $\bm{U}\in W^{1,2}(B_R, \rNn)$, and a sequence $\{\ep _k\}$ such that $\ep_k \to 0^+$,
\begin{equation}\label{loc19}
b_{\ep _k}(|\nabla \bu_{\ep _k}|)^{p-2} \nabla \bu_{\ep _k} \to \bm{U}\quad \hbox{in $L^2(B_R, \rNn)$} \quad \hbox{and} \quad b_{\ep _k}(|\nabla \bu_{\ep _k}|)^{p-2} \nabla \bu_{\ep _k} \rightharpoonup \bm{U}\quad \hbox{in $W^{1,2}(B_R, \rNn)$}
\end{equation}
as $k \to \infty$.
Since $\lim _{k\to \infty} b_{\ep _k}(t)^{p-2} = t^{p-2}$ for $t>0$,
thanks to \eqref{loc19} and \eqref{akm4},
\begin{equation}\label{loc22}
|\nabla \bm{u}|^{p-2} \nabla \bm{u} = \bm{U} \in W^{1,2}(B_R, \rNn).
\end{equation}
Moreover, equations \eqref{loc18} and  \eqref{loc19} entail that
\begin{align}\label{loc26}
\| |\nabla \bu|^{p-2}\nabla \bu\|_{W^{1,2}(B_R, \rNn)}   \leq C\big(\|\bm{f}\|_{L^2(B_{2R}, \rN)} + (R^{-\frac n2}+R^{-\frac n2-1})\|\nabla \bu\|_{L^{p-1}(B_{2R}, \rNn)}^{p-1}\big)\,.
\end{align}
It remains to remove assumption \eqref{fsmoothloc}. Suppose that
 $\bm{f}\in L^2_{\rm loc}(\o, \rN)$, let $\bu$ be an approximable local solution  to equation \eqref{localeq}, and let $\bm{f}_k$ and $\bu_k$ be as in the  definition  of this kind of solution given at the begining of the present section.
An application of inequality \eqref{loc26}  to $\bu_k$ tells us that 
$|\nabla \bu_k|^{p-2} \nabla \bu_k   \in W^{1,2}(B_R, \rNn)$, and 
\begin{align}\label{loc27}
\||\nabla \bu_k|^{p-2}\nabla \bu_k\|_{W^{1,2}(B_R, \rNn)}&  \leq C\big(\|\bm{f}_k\|_{L^2(B_{2R}, \rN)} + (R^{-\frac n2}+ R^{-\frac n2-1})\|\nabla \bu_k\|_{L^{p-1}(B_{2R}, \rNn)}^{p-1}\big)
\,,
\end{align}
where the constant $C$ is independent of $k$. Therefore, owing to equation \eqref{approxaloc}, the sequence $\{|\nabla \bu_k|^{p-2}\nabla \bu_k\}$ is bounded in $W^{1,2}(B_R, \rNn)$, and hence 
there exists a function
$\bm{U} \in W^{1,2}(B_R, \rNn)$, and a subsequence, still indexed by $k$,   such that
\begin{equation}\label{loc28}
|\nabla \bu_k|^{p-2} \nabla \bu_k \to \bm{U}\quad \hbox{in $L^2(B_R, \rNn)$} \quad \hbox{and} \quad |\nabla \bu_k|^{p-2}\nabla \bu_k \rightharpoonup \bm{U}\quad \hbox{in $W^{1,2}(B_R, \rNn)$}.
\end{equation}
By assumption \eqref{approxuloc}, $\nabla \bu_k \to \nabla \bu$ a.e. in $\Omega$. Hence, owing to \eqref{loc28}, 
\begin{equation}\label{loc30}
|\nabla \bu|^{p-2} \nabla \bu = \bm{U}\in W^{1,2}(B_R, \rNn)\,.
\end{equation}
Inequality \eqref{secondloc2} follows from \eqref{loc27}, via \eqref{approxaloc}, \eqref{loc28}  and \eqref{loc30}. \qed

\section{Global estimates}\label{proofs}

This section is devoted to the proofs of   Theorems  \ref{seconddir} and \ref{secondconvex} and Corollary \ref{secondmarc}. As a preliminary, we recall some relevant notions of   solutions to Dirichlet and Neumann boundary value problems. 
\par We begin with Dirichlet problems. Let $\o$ be an open set with finite Lebesgue measure $|\o|$. Assume that $\bm f \in L^1(\Omega, \rN) \cap (W^{1,p}_0(\Omega, \rN))'$. A function $\bu \in W^{1,p}_0(\o, \mathbb R^N)$ is called a weak solution to the Dirichlet problem \eqref{localeq}$+$\eqref{dircond} if  
\begin{equation}\label{weaksol}
\int _\o |\nabla \bu|^{p-2} \nabla \bu \cdot \nabla  \bm{\varphi} \color{black} \,dx = \int _{\Omega} {\bf f}\cdot \bm{\varphi}  \,dx
\end{equation}
for every $\bm{\varphi}  \in W^{1,p}_0(\Omega, \mathbb R^N)$. Classically, a unique weak solution to \eqref{localeq}$+$\eqref{dircond}  exists under the present assumptions on $\o$ and  $\bm f$.
\\ 
Assume next that $\bm f \in L^q(\Omega, \rN)$ for some $q \geq 1$. 
A function $\bu \in W^{1,1}_0(\o, \mathbb R^N)$ is called an approximable solution to the Dirichlet problem \eqref{localeq}$+$\eqref{dircond} if   there exists a sequence $\{{\bf f}_k\} \subset C^{\infty}_0(\Omega, \mathbb R^N)$ such that  ${\bf f}_k \to \bm{f}$ in $L^q(\o, \mathbb R^N)$, and 
the sequence $\{\bu _k\}$ of weak solutions to the Dirichlet problems 
\begin{equation}\label{eqdirichletk}
\begin{cases}
-{\rm {\bf  div}} (|\nabla \bu _k|^{p-2}\nabla {\bf u} _k) = {\bf f}_k  & {\rm in}\,\,\, \o \\
 {\bf u}_k =0  &
{\rm on}\,\,\,
\partial \o \,
\end{cases}
\end{equation}
satisfies
\begin{equation}\label{convdir}
\bu _k \to \bu \quad \hbox{and} \quad \nabla \bu _k \to \nabla \bu \quad \hbox{a.e. in $\o$.}
\end{equation}

\smallskip
\par
Parallel definitions can be given for Neumann problems. Let $\o$ be a bounded Lipschitz domain. Assume that $\bm f \in L^1(\Omega, \rN) \cap (W^{1,p}(\Omega, \rN))'$ and satisfies the compatibility condition \eqref{0mean}. 
A function $\bu \in W^{1,p}(\o, \mathbb R^N)$ is called a weak solution to the Neumann problem \eqref{localeq}$+$\eqref{neucond}  if  
\begin{equation}\label{weaksoneul}
\int _\o |\nabla \bu|^{p-2} \nabla \bu \cdot \nabla  \bm{\varphi} \color{black} \,dx = \int _{\Omega} {\bf f}\cdot \bm{\varphi}  \,dx
\end{equation}
for every $\bm{\varphi}  \in W^{1,p}(\Omega, \mathbb R^N)$. A weak solution to  \eqref{localeq}$+$\eqref{neucond} exists if $\o$ and $\bm f$ iare as above, and is unique, up to additive constant vectors in $\rN$.
\\ Assume now that $\bm f \in L^q(\Omega, \rN)$ for some $q \geq 1$. 
A function $\bu \in W^{1,1}(\o, \mathbb R^N)$ is called an approximable solution to the Neumann problem  \eqref{localeq}$+$\eqref{neucond}  if  there exists a sequence $\{{\bf f}_k\} \subset C^{\infty}_0(\Omega, \mathbb R^N)$ such that $\int _\o \bm f_k\,dx =0$ for $k \in \mathbb N$, ${\bf f}_k \to \bm{f}$ in $L^q(\o, \mathbb R^N)$, and 
a sequence $\{\bu _k\}$ of weak solutions to the Neumann problems 
\begin{equation}\label{eqneumannk}
\begin{cases}
- {\rm {\bf  div}}(|\nabla \bu _k|^{p-2}\nabla {\bf u} _k) = {\bf f}_k & {\rm in}\,\,\, \o \\
  \displaystyle  \frac{\partial{\bf u}_k}{\partial \bm \nu} =0  &
{\rm on}\,\,\,
\partial \o \,
\end{cases}
\end{equation}
satisfies
\begin{equation}\label{convneu}
\bu _k \to \bu \quad \hbox{and} \quad \nabla \bu _k \to \nabla \bu \quad \hbox{a.e. in $\o$.}
\end{equation}

The following a priori estimate for the gradient of weak solutions to Dirichlet and Neumann problems will be of use in the proof of our global main results.

\begin{proposition}\label{talentil2}
 Assume that  $n\geq 2$, $N \geq 1$ and $p >1$. Let $\o$ be an open set in $\rn$ such that $|\o|<\infty$.
\\ 
{\rm (i)}  Assume that $\bm{f} \in L^q(\o, \rN)\cap (W^{1,p}_0(\o, \rN))'$ for some $q>1$. Let $\bu$ be the weak solution to the Dirichlet problem \eqref{localeq}$+$\eqref{dircond}. Then, there exists a constant $C=C(n, N, p,q, |\o|)$ such that
\begin{equation}\label{talentigrad}
\| \nabla \bu\|_{L^{p-1}(\o , \rNn)} \leq C\|\bm{f}\|_{L^q(\o, \rN)}^{\frac 1{p-1}}\,.
%
%
\end{equation}
\\ 
{\rm (ii)}  Suppose, in addition, that $\Omega$ is a bounded  Lipschitz domain. Assume that $\bm{f} \in L^q(\o, \rN)\cap (W^{1,p}(\o, \rN))'$ and satifies condition \eqref{0mean}. Let $\bu$ be a weak solution to the Neumann problem \eqref{localeq}$+$\eqref{neucond}. Then there exists a constant $C=C(n,N, p,q,L_\o,d_\o)$ such that
\begin{equation}\label{neumanngrad}
\| \nabla \bu\|_{L^{p-1}(\o , \rNn)} \leq C\|\bm{f}\|_{L^q(\o, \rN)}^{\frac 1{p-1}}\,.
\end{equation}
\end{proposition}

\par\noindent
{\bf Proof}. Part (i). Given $\delta >0$ and $\gamma \in (0,1)$, choose the test function $(|\bu|+\delta)^{-\gamma} \bu$ in the definition of weak solution to problem  \eqref{localeq}$+$\eqref{dircond}.
Since, by the chain rule for Sobolev functions,
\begin{equation}\label{talentil21}
\big((|\bu|+\delta)^{-\gamma} \bu)_{x_i} = (|\bu|+\delta)^{-\gamma} \bu_{x_i} -\gamma  (|\bu|+\delta)^{-\gamma-1}\frac{\bu}{|\bu|}\bu \cdot \bu_{x_i} \quad \hbox{a.e. in $\o$,}
\end{equation}
for $i=1, \dots , n$, 
we obtain that
\begin{align}\label{talentil22}
\int _\o & |\nabla \bu|^{p-2} \sum _{i=1}^n \bigg[(|\bu|+\delta)^{-\gamma} \bu_{x_i} -\gamma  (|\bu|+\delta)^{-\gamma-1}\frac{\bu}{|\bu|}\bu \cdot \bu_{x_i}\bigg]\cdot \bu_{x_i}\, dx 
\\ \nonumber & = 
\int _\o |\nabla \bu|^{p-2}  \sum _{i=1}^n \bigg[(|\bu|+\delta)^{-\gamma} |\bu_{x_i}|^2 -\gamma  (|\bu|+\delta)^{-\gamma-1}\frac{(\bu \cdot \bu_{x_i})^2}{|\bu|}\bigg]\, dx
\\ \nonumber &  = \int _\o (|\bu|+\delta)^{-\gamma} {\bf f} \cdot \bu \, dx\,.
%
%
\end{align}
Hence,
\begin{align}\label{talentil23}
\int _\o   |\nabla \bu|^p \big[(|\bu|+\delta)^{-\gamma}   -\gamma  (|\bu|+\delta)^{-\gamma-1}|\bu|\big]\, dx
\leq \int _\o (|\bu|+\delta)^{-\gamma} |{\bf f}| |\bu| \, dx\,.
\end{align}
Passing to the limit as $\delta \to 0^+$ in \eqref{talentil23} yields, via monotone convergence, 
\begin{align}\label{talentil24}
(1-\gamma)\int _\o    |\nabla \bu|^p |\bu|^{-\gamma}   \, dx
\leq \int _\o |{\bf f}| |\bu|^{1-\gamma}  \, dx\,.
\end{align}
Notice that the right-hand side of inequality \eqref{talentil24} is finite, if $\gamma$ is sufficiently close to $1$, by H\"older's inequality and the Sobolev embedding theorem, since $\bm{f}\in L^q(\o, \rN)$ and $\bu \in W^{1.p}_0(\o, \rN)$. As a consequence of \eqref{talentil24}, one can thus show that $|\bu|^{-\frac \gamma p}\bu \in W^{1.p}_0(\o, \rN)$. To verify this assertion, observe that
\begin{align}\label{talentil35}
\big|\nabla 
\big((|\bu|+\delta)^{-\frac \gamma p}\bu\big)\big|^p \leq C |\nabla \bu|^p (|\bu |+\delta)^{-\gamma} \leq |\nabla \bu|^p |\bu|^{-\gamma} \quad \hbox{a.e. in $\o$,}
\end{align}
for every $\delta >0$, and for some constant $C=C(n, N, p, \gamma)$. Therefore, by \eqref{talentil24} the family of functions $\{(|\bu|+\delta)^{-\frac \gamma p}\bu\}$ is bounded in $W^{1.p}_0(\o, \rN)$ as $\delta \to 0^+$. Hence, by weak compactness, its   limit  $|\bu|^{-\frac \gamma p}\bu \in W^{1.p}_0(\o, \rN)$.
\\ Let $r \in (p, \tfrac {np}{n-p}]$ if $p<n$, or $r>p$ if $p \geq n$.
The Poincar\'e-Sobolev inequality applied to the function $|\bu|^{-\frac \gamma p}\bu$ tells us that
\begin{align}\label{talentil28}
\bigg(\int _\o |\bu|^{(1-\frac \gamma p)r}\, dx\bigg)^{\frac 1r} = \bigg(\int _\o \big||\bu|^{-\frac \gamma p}\bu|^r\, dx\bigg)^{\frac 1r} & \leq  C \bigg(\int _\o \big|\nabla 
\big(|\bu|^{-\frac \gamma p}\bu\big)\big|^p\, dx\bigg)^{\frac 1p}
 \\ \nonumber & \leq C'
\bigg(\int _\o|\nabla \bu|^p |\bu|^{-\gamma}\, dx\bigg)^{\frac 1p}\,
\end{align} 
for some constants  $C=C(p,n,r, |\o|)$ and  $C'=C'(p,n,r, \gamma, |\o|)$.
If 
\begin{equation}\label{cond1}
q'(1-\gamma) \leq (1-\tfrac \gamma p)r\,,
\end{equation}
then
from  inequality \eqref{talentil24}, H\"older's inequality and inequality \eqref{talentil28} one can infer that
\begin{align}\label{talentil29}
(1-\gamma) \int _\o  |\nabla \bu|^p |\bu|^{-\gamma}   \, dx 
 & \leq \bigg( \int _\o |{\bf f}|^q\, dx\bigg)^{\frac 1q} \bigg(\int_\o |\bu|^{q'(1- \gamma)}\bigg)^{\frac{1}{q'}}
\\ \nonumber & \leq C \bigg( \int _\o |{\bf f}|^q\, dx\bigg)^{\frac 1q} \bigg(\int_\o |\bu|^{(1-\frac \gamma p)r}\bigg)^{\frac{1-\gamma}{(1-\frac{\gamma}p)r}}
\\ & \nonumber  \leq C'   \bigg( \int _\o |{\bf f}|^q\, dx\bigg)^{\frac 1q} \bigg(\int_\o|\nabla \bu|^p |\bu|^{-\gamma}\, dx\bigg)^{\frac{1-\gamma}{p-\gamma}}\,,
\end{align}
for some constants $C=C(n,N,p,q,r,\gamma,|\o|)$ and $C'=C'(n,N,p,q,r,\gamma,|\o|)$.
Hence,
\begin{align}\label{talentil30}
 \bigg(\int _\o  |\nabla \bu|^p |\bu|^{-\gamma}   \, dx \bigg)^{\frac{p-1}{p-\gamma}}
  \leq C   \bigg( \int _\o |{\bf f}|^q\, dx\bigg)^{\frac 1q}\,
\end{align}
for some constant $C=C(n,N,p,q,r,\gamma,|\o|)$.
On the other hand, if
\begin{equation}\label{cond2}
\gamma (p-1) \leq (1-\tfrac \gamma p)r\,,
\end{equation}
 then the use of H\"older's inequality (twice), and  of inequality \eqref{talentil28} yields
\begin{align}\label{talentil32}
\int _\o |\nabla \bu|^{p-1}\, dx & \leq 
 \bigg(\int _\o  |\nabla \bu|^p |\bu|^{-\gamma}   \, dx \bigg)^{\frac{1}{p'}} \bigg(\int _\o   |\bu|^{\gamma(p-1)}   \, dx \bigg)^{\frac{1}{p}}
\\ \nonumber & 
  \leq C    \bigg(\int _\o  |\nabla \bu|^p |\bu|^{-\gamma}   \, dx \bigg)^{\frac{1}{p'}} \bigg(\int _\o |\bu|^{(1-\frac \gamma p)r}\, dx\bigg)^{\frac{\gamma (p-1)}{(p-\gamma)r}}
\\ \nonumber & 
  \leq C'    \bigg(\int _\o  |\nabla \bu|^p |\bu|^{-\gamma}   \, dx \bigg)^{\frac{p-1}{p-\gamma}} \,
\end{align}
for some constants $C=C(n,N,p,r,\gamma,|\o|)$ and $C'=C'(n,N,p,r,\gamma,|\o|)$. 
Notice that both conditions \eqref{cond1} and \eqref{cond2} are fulflilled owing to our choice of $r$, provided that $\gamma$ is sufficiently close to $1$. 
Inequality \eqref{talentigrad} follows from \eqref{talentil30} and \eqref{talentil32}.
\\ Part (ii). The proof follows along the same lines as that of Part (i). One has just to replace  (if necessary) the original solution $\bu$ to the Neumann problem \eqref{localeq}$+$\eqref{neucond},  by the solution $\bu- \bm\zeta_0$,  where $\bm\zeta_0 \in \rN$  is chosen in such  a way that
\begin{equation}\label{normal}
\int _\o |\bu - \bm\zeta_0 |^{-\frac \gamma p} (\bu - \bm\zeta_0)\, dx = 0\,.
\end{equation}
\color{black}
To verify  that this choice is indeed possible, consider the function $\Phi  : \rN \to [0, \infty)$, given by 
\begin{equation}\label{Phi}
\Phi (\bm\zeta) = \int _\o  |\bu- \bm\zeta|^{2-\frac \gamma p} \, dx \quad \hbox{for $\bm\zeta \in \rN$.}
\end{equation}
This function is actually finite-valued if $\gamma$ is sufficiently close to $1$, since $\bm u \in W^{1,p}(\o, \rN)$, and hence $\bm u \in L^{2-\frac \gamma p}(\o, \rN)$ by the Sobolev embedding theorem. Moreover, $\Phi$  is (strictly) convex, continuosly differentiable, and $\lim _{|\bm\zeta|\to \infty}\Phi (\bm\zeta) =\infty$.
%
%
Thus, there exists a (unique) minimum point $\bm\zeta_0$ of $\Phi$, and 
\begin{equation}\label{critical}
\big(2-\tfrac \gamma p\big)\int _\o |\bu- \bm\zeta_0|^{-\frac \gamma p} (\bu - \bm\zeta_0)\, dx = \nabla \Phi (\bm\zeta_0) =0\,.
\end{equation}
On denoting the solution $\bu - \bm\zeta_0$ simply by $\bu$, condition \eqref{normal} reads 
\begin{equation}\label{normalu}
\int _\o |\bu |^{-\frac \gamma p} \bu\, dx = 0\,.
\end{equation}
Under condition \eqref{normalu}, inequality \eqref{talentil28} continues to hold, by the Poincar\'e-Sobolev inequality for functions with vanishing mean value, for some constants $C$ and $C'$. The dependence of these constants   on $\o$ is as described in the last part of the statement. Having disposed of inequality \eqref{talentil28}, inequality \eqref{neumanngrad} can be derived via the same argument as in the proof of \eqref{talentigrad}.
\qed

Our proof of Theorem \ref{seconddir} entails an approximation of the domain $\o$ by a sequence of domains $\o_m$ with a smooth boundary. Such an approximation requires, in particular, that the quantities $L_{\o_m}$ and $d_{\o_m}$, and the functions $\mathcal K_{\o_m}$ be uniformly bounded (up to multiplicative constants) by $L_{\o}$ and $d_{\o}$, and  $\mathcal K_{\o}$. This is the subject of the next lemma. The main steps of its proof are described in the Appendix.

\begin{lemma}\label{approxcap}
Let $\Omega$ be a bounded Lipschitz domain in $\rn$, $n \geq 2$ such that $\partial \o \in W^{2,1}$. Assume that the function $\mathcal K_\o (r)$, defined as in \eqref{defK}, is finite-valued for $r\in (0,1)$.
Then there exist positive constants $r_0$ and $C$ and a sequence of bounded open sets $\{\Omega_m\}$, 
such that $\partial \Omega _m \in C^\infty$, $\Omega \subset \Omega _m$, $\lim _{k \to \infty}|\Omega _m \setminus \Omega| = 0$, the Hausdorff distance between $\Omega _m$ and $\Omega$ tends to $0$ as $m \to \infty$,
\begin{equation}\label{feb35}
L_{\o _m} \leq C L_\o \,, \quad d_{\o _m} \leq C d_\o
\end{equation}
and
\begin{equation}\label{appcap0}
\mathcal K_{\o_m}(r) \leq C \mathcal K_{\o} (r)
\end{equation}
for $r\in (0, r_0)$ and $m \in \mathbb N$.
\end{lemma}

\smallskip
\noindent
{\bf Proof of Theorem \ref{seconddir}}. We begin by dealing with the case when $\bu$ is a solution to the Dirichlet problem \eqref{localeq}$+$\eqref{dircond}. The proof relies on a combination of Theorem \ref{step1}, Part (ii), with   approximation arguments for  the differential operator, the domain and the right-hand side. The  relevant approximations  are accomplished in separate steps.
\par
\noindent \emph{Step 1}. Here, we assume 
the  additional conditions 
 \begin{equation}\label{fsmooth}
\bm{f} \in C^\infty_0(\o, \rN)\,,
\end{equation}
 and 
 \begin{equation}\label{omegasmooth}
\partial \Omega \in C^\infty\,.
\end{equation}
Given $\ep \in (0,1)$,  denote by $\bu_\varepsilon$  the weak   solution to the system
\begin{equation}\label{eqdirichletep}
\begin{cases}
- {\rm {\bf  div}}(b_\varepsilon  (|\nabla \bu_\varepsilon|)^{p-2} \nabla {\bf u}_\varepsilon ) = {\bf f} & {\rm in}\,\,\, \o \\
 {\bf u_\varepsilon} =0  &
{\rm on}\,\,\,
\partial \o \,,
\end{cases}
\end{equation}
with $b_\ep$ defined as in \eqref{aeps}.
By \cite[Theorem 2.1]{cmARMA}, there exists a constant $C$, independent of $\ep$, such that 
\begin{align}\label{gradbound}
\|\nabla \bu_\varepsilon\|_{L^\infty (\o, \rNn)} \leq C.
\end{align}
 Hence, for each $\varepsilon \in (0,1)$, there exist constants $c_2>c_1>0$ such that
\color{black}
\begin{equation}\label{abound}
 c_1 \leq b_\varepsilon (|\nabla \bu_\varepsilon|) \leq c_2 \qquad \hbox{in $\o$.}
\end{equation}
\\ Properties \eqref{fsmooth}, \eqref{omegasmooth} and \eqref{abound} permit an application of  a result by Elcrat and Meyers \cite[Theorem 8.2]{BF},  ensuring that ${\bf u_\varepsilon}\in W^{2,2}(\Omega, \rN )$. 
Thus, ${\bf u}_\varepsilon \in W^{1,2}_0(\Omega, \rn) \cap W^{1,\infty}(\Omega, \rN) \cap
W^{2,2}(\Omega, \rN)$. By  standard approximation \cite[Chapter 2, Corollary 3]{Burenkov},  there exists a
sequence $\{{\bf u}_k\} \subset C^\infty (\Omega, \rN)\cap
C^2(\overline \Omega, \rN)$ such that $\bu_k = 0$ on $\partial
\Omega$ for $k \in \mathbb N$, and
\begin{equation}\label{convk}
{\bf u}_k \to {\bf u}_\varepsilon \quad \hbox{in $W^{1,2}_0(\Omega, \rN)$,}
\quad {\bf u}_k \to {\bf u}_\varepsilon \quad \hbox{in $W^{2,2}(\Omega, \rN)$,}
\quad \nabla {\bf  u}_k \to \nabla {\bf u}_\varepsilon 
\quad \hbox{a.e. in
$\Omega $},
\end{equation}
as $k \to \infty$. Furthermore, 
\begin{equation}\label{boundk}
\|\nabla \bu_k\|_{L^\infty (\o, \rNn)}  \leq C \|\nabla \bu_\varepsilon\|_{L^\infty (\o, \rNn)} 
\end{equation}
\color{black}
for some constant $C$ independent of $k$, 
and,
%
%
 by the chain rule for
vector-valued Sobolev functions \cite[Theorem 2.1]{MarcusMizel}, $|\nabla |\nabla {\bf u}_k|| \leq |\nabla ^2
\bu_k|$ a.e. in $\Omega$.
\\ Finally, one can  show that
\begin{align}\label{sys1k}
- {\rm {\bf div}} (b_\varepsilon(|\nabla {\bf u_k}|)^{p-2}\nabla {\bf u_k} ) \to {\bf
f}
\quad \hbox{in $L^2(\Omega , \rN)$},
\end{align}
as $k \to \infty$, see  \cite[Equation (6.12)]{cmARMA}. Owing to assumption \eqref{capcond}, inequality \eqref{fund}, applied with $a$ replaced by $b_\varepsilon ^{p-2}$ and $\bu$ replaced by $\bu _k$, yields
\begin{multline}\label{main16k}
\| b_\varepsilon(|\nabla \bu_k|)^{p-2}\nabla \bu_k\|_{W^{1,2}(\o, \rNn)} \\ \leq C \Big(
\|{\rm  {\bf div}} ( b_\varepsilon(|\nabla \bu_k|)^{p-2} \nabla \bu_k)\|_{L^2(\o, \rN)} +
 \| b_\varepsilon(|\nabla \bu_k|)^{p-2}\nabla \bu_k\|_{L^1(\o, \rNn)}\Big)
\end{multline}
for $k\in \mathbb N$, and for some constant $C=C(n, N, p, L_\o, d_\o, \mathcal K_\o)$.  Observe that this constant is independent of $\ep$, thanks to \eqref{epsindex}.
Owing to equations  \eqref{boundk}--\eqref{main16k}, the sequence $\{b_\varepsilon(|\nabla \bu_k|)^{p-2}\nabla \bu_k\}$ is bounded in $W^{1,2}(\o, \rNn)$. Hence, there exists a subsequence of $\{\bu_k\}$, still denoted by $\{\bu_k\}$, and a function ${\bf U}_\varepsilon \in W^{1,2}(\o, \rNn)$ such that
\begin{equation}\label{main32keps}
b_\varepsilon(|\nabla \bu_{k}|)^{p-2} \nabla \bu_{k} \to {\bf U}_\varepsilon\quad \hbox{in $L^2(\Omega, \rNn )$,} \quad b_\varepsilon(|\nabla \bu_{k}|)^{p-2} \nabla \bu_{k} \rightharpoonup {\bf U}_\varepsilon\quad \hbox{in $W^{1,2}(\Omega, \rNn )$}.
\end{equation}
Since, by \eqref{convk}, $\nabla \bu_k \to \nabla \bu_\varepsilon$ a.e. in $\o$, one has that
\begin{equation}\label{main31keps}
b_\varepsilon(|\nabla \bu_{k}|)^{p-2} \nabla \bu_{k} \to b_\varepsilon(|\nabla \bu_\varepsilon|) ^{p-2}\nabla \bu_\varepsilon \quad \hbox{a.e. in $\Omega$.}
\end{equation}
Coupling \eqref{main31keps} with \eqref{main32keps} ensures that
\begin{equation}\label{main100eps}
b_\varepsilon(|\nabla \bu_\varepsilon|)^{p-2} \nabla \bu_\varepsilon  = {\bf U}_\varepsilon \in W^{1,2}(\Omega, \rNn )\,,
\end{equation}
and passing to the limit as $k \to \infty$ in \eqref{main16k}
yields
\begin{align}\label{main16eps}
\| b_\varepsilon(|\nabla \bu_\varepsilon|)^{p-2}\nabla \bu_\varepsilon\|_{W^{1,2}(\o, \rNn)} \leq C \big(
\|{\bf f}\|_{L^2(\o, \rN) } + 
 \|b_\varepsilon(|\nabla \bu_\varepsilon|)^{p-2}\nabla \bu_\varepsilon\|_{L^1(\o, \rNn)}\big)\,.
\end{align}
Note that, in deriving inequality \eqref{main16eps}, we have made use of \eqref{main32keps} and \eqref{main100eps} on the left-hand side, and of   \eqref{boundk} and \eqref{sys1k} on the right-hand side.  From inequalities \eqref{main16eps} and \eqref{gradbound}, we infer that there exists a constant $C$, independent of $\ep$,   such that 
\begin{align}\label{main16eps'}
\| b_\varepsilon(|\nabla \bu_\varepsilon|)^{p-2}\nabla \bu_\varepsilon\|_{W^{1,2}(\o, \rNn)} \leq C (1 + \ep ^{p-1})\,.
\end{align}
 Owing to inequality \eqref{main16eps'}, the family of functions  $b_\varepsilon(|\nabla \bu_\varepsilon|)^{p-2}\nabla \bu_\varepsilon$ is uniformly bounded in $W^{1,2}(\o, \rNn)$ for $ \ep \in (0,1)$, and hence there exists a sequence $\{\varepsilon _m\}$ and a function ${\bf U}  \in W^{1,2}(\o, \rNn)$ such that $\ep _m \to 0$, and
\begin{equation}\label{main32eps}
b_{\varepsilon _m}(|\nabla \bu_{\varepsilon _m}|)^{p-2} \nabla \bu_{\varepsilon _m} \to {\bf U} \,\,\hbox{in $L^2(\Omega, \rNn)$,} \quad b_{\varepsilon _m}(|\nabla \bu_{\varepsilon _m}|)^{p-2} \nabla \bu_{\varepsilon _m} \rightharpoonup {\bf U} \, \,\hbox{in $W^{1,2}(\Omega, \rNn )$}.
\end{equation}
Now, an analogous (and even simpler) argument as in the proof of \eqref{akm4} tells us that
\begin{align}\label{akmglob}
\nabla \bu _{\ep_m} \to \nabla \bu \qquad \hbox{in $L^p(\o, \rNn)$.}
\end{align}
In particular, notice that, in this argument, inequality \eqref{cialpha1} has to be replaced by
\begin{align*}
\int _{\o} |\nabla \bu _{\ep_m}|^p\, dx \leq C \bigg(\int _{\o}|\bm{f}|^{p'}\,dx + \ep_m^p\bigg)\,,
\end{align*}
an easy consequence of the use of $\bu _{\ep_m}$ as a test function in the definition of weak solution to problem \eqref{eqdirichletep}, with $\ep = \ep_m$.
From equations \eqref{main32eps} and \eqref{akmglob} one infers that
\begin{equation}\label{main25}
|\nabla \bu|^{p-2} \nabla \bu = {\bf U} \in W^{1,2}(\Omega, \rNn)\,.
\end{equation}
Also, equation \eqref{akmglob}, inequality \eqref{gradbound},  the dominated convergence theorem for Lebesgue integrals and 
inequality \eqref{talentigrad} ensure that
\begin{align}\label{napoli}
\lim _{m \to \infty}  \|b_{\varepsilon _m}(|\nabla \bu_{\varepsilon_m}|)^{p-2}\nabla \bu_{\varepsilon_m}\|_{L^1(\o, \rNn)} = 
\|\nabla \bu\|_{L^{p-1}(\o, \rNn)} \leq C \|\bm{f}\|_{L^2(\o, \rN)}
\end{align}
for some constant $C=C(n,N,p, |\o|)$. 
From \eqref{main16eps}, we obtain via \eqref{main32eps}, \eqref{main25} and \eqref{napoli} that
\begin{equation}\label{main27}
\||\nabla \bu|^{p-2} \nabla \bu\|_{W^{1,2}(\o, \rNn)} \leq C \|{\bf f}\|_{L^2(\o, \rN)}
\end{equation}
%
%
for some constant $C= C(n, N, p, L_\o, d_\o, \mathcal K_\o)$.

\smallskip
\par\noindent
\emph{Step 2}. Here, we remove assumption \eqref{omegasmooth}.  Let $\{\o_m\}$ be a sequence of open sets approximating $\o$  in the sense of Lemma \ref{approxcap}. Consider, for each $m \in \mathbb N$, the weak solution  $\bm u_m$  to the Dirichlet problem 
\begin{equation}\label{eqm}
\begin{cases}
- {\rm {\bf div}} (|\nabla \bu_m|^{p-2} \nabla \bu_m ) = {\bf f}  & {\rm in}\,\,\, \o _m \\
 \bu_m =0  &
{\rm on}\,\,\,
\partial \o _m \,,
\end{cases}
\end{equation}
where ${\bf f}$ still fulfils \eqref{fsmooth}, and  is extended by $0$ outside $\Omega$. By inequality \eqref{main27} of Step 2,  applied to $\bm u_m$,
\begin{align}\label{main27aus}
\||\nabla \bu_m|^{p-2} \nabla \bu_m\|_{W^{1,2}(\o, \rNn)} & \leq 
\||\nabla \bu_m|^{p-2} \nabla \bu_m\|_{W^{1,2}(\o_m, \rNn)} \\ \nonumber & \leq C \|{\bf f}\|_{L^2(\o_m, \rN)}= C \|{\bf f}\|_{L^2(\o, \rN)},
\end{align}
%
%
for some constant  $C(n, N, p, L_\o, d_\o, \mathcal K_\o)$. 
Note that this dependence of the constant $C$ is guaranteed by 
properties  \eqref{feb35} and \eqref{appcap0} of the sequence $\{\o_m\}$.
\\
Thanks to \eqref{main27aus}, the sequence  
$\{|\nabla \bu_m |^{p-2} \nabla  \bu_m\}$ is  bounded in $W^{1,2}(\Omega, \rNn)$,
and hence  there exists a subsequence, still denoted by $\{\bu_m\}$ and a function  ${\bf U} \in W^{1,2}(\Omega, \rNn)$, such that
\begin{equation}\label{main32}
|\nabla \bu_{m}|^{p-2} \nabla \bu_{m} \to {\bf U}\quad \hbox{in $L^2(\Omega, \rNn)$,} \qquad   |\nabla \bu_{m}|^{p-2} \nabla \bu_{m} \rightharpoonup {\bf U}\quad \hbox{in $W^{1,2}(\Omega, \rNn)$}.
\end{equation}
Next, we claim that 
there exists $\alpha \in (0,1)$ such that $\bu _m\in C^{1,\alpha}_{\rm loc}(\Omega, \rN)$, and for every open set $\o ' \subset \subset \Omega$ there exists a constant $C$, independent of $m$, such that  
\begin{align}\label{main29}
\|\bu_{m}\|_{C^{1,\alpha}(\Omega ', \rN)} \leq C\,.
\end{align}
 Actually, it follows 
from a special case of   \cite[Corollary 5.6]{BCDKS} that, for each open set $\o'$ as above, there exists a constant $C$, independent of $m$, such that
\begin{align}\label{Liploc}
\|\nabla \bu _{m}\|_{L^\infty(\Omega ', \rNn)} \leq C.
\end{align}
 Thanks to \cite[Corollary 1.26 and equation (1.8d)]{sebastian} and inequality \eqref{Liploc},  
\begin{align}\label{calpham}
\|\nabla \bu _{m}\|_{C^{\alpha}(\Omega ', \rNn)} \leq C\,
\end{align}
for some constant $C$ independent of $m$. 
Since assumtpion \eqref{fsmooth} is still in force, a basic energy estimate yields
\begin{align}\label{energyst}
\|\nabla \bu _{m}\|_{L^p(\Omega_m , \rNn)} \leq C
\end{align} 
for some constant $C$ independent of $m$. Hence, by the Poincar\'e  inequality, 
\begin{align}\label{energym}
 \| \bu _{m}\|_{L^p(\Omega_m , \rN)} \leq C\,,
\end{align}
 with $C$ indepedent of $m$, inasmuch as $\bu _{m} \in W^{1,p}_0(\o_m, \rN)$, and $\o _m$ satisfies \eqref{feb35}. 
%
%
Owing to \eqref{Liploc} and \eqref{energym},   a Sobolev type inequality yields
\begin{align}\label{linfinitym}
 \| \bu _{m}\|_{L^\infty(\Omega ', \rN)} \leq C\,
\end{align}
for some constant $C$ independent of $m$.
In this connection, notice  that, without loss of generality, $\partial \o'$ can be assumed to be smooth. 
Inequality \eqref{main29} follows from \eqref{calpham} and \eqref{linfinitym}.
\color{black}
\\ Thus,  on taking, if necessary,  a further subsequence,  
\begin{equation}\label{main30}
\bu_{m} \to {\bf v}\, \quad \hbox{and}\quad \nabla \bu_{m} \to \nabla {\bf v} \quad \hbox{in $\Omega$,}
\end{equation}
for some function ${\bf v} \in C^1(\o, \rN)$. 
In particular,
\begin{equation}\label{main31}
|\nabla \bu_{m}|^{p-2} \nabla \bu_{m} \to |\nabla {\bf v}|^{p-2} \nabla {\bf v} \quad \hbox{in $\Omega$.}
\end{equation}
By \eqref{main31} and \eqref{main32}, 
\begin{equation}\label{main100}
|\nabla {\bf v}|^{p-2} \nabla {\bf v}   = {\bf U} \in W^{1,2}(\Omega, \rNn)\,.
\end{equation}
 Now,  let us make use of a test  function $\bm{\varphi} \in C^\infty _0(\o, \rN)$ (extended by $0$ outside $\o$) in the  definition of weak solution to problem \eqref{eqm}, and 
  pass to the limit as $m \to \infty$ in the resulting equation, namely in the equation
\begin{equation}\label{main101}
\int _{\Omega_m} |\nabla \bu_m|^{p-2} \nabla \bu_m \cdot \nabla \bm{\varphi} \, dx = \int_{\Omega _m} {\bf f} \cdot \bm{\varphi} \, dx\,.
\end{equation}
So doing, we infer from \eqref{main32} and \eqref{main100} that
\begin{equation}\label{main101bis}
\int _\Omega |\nabla {\bf v}|^{p-2} \nabla {\bf v} \cdot \nabla \bm{\varphi} \, dx = \int_\Omega {\bf f} \cdot \bm{\varphi} \, dx\,.
\end{equation}
Inequality \eqref{energyst} ensures that $\|\nabla \bu _m\|_{L^p(\o, \rNn)} \leq C$ for some constant $C$ independent of $m$. The same inequality thus holds with $\bu _m$ replaced with ${\bf v}$, 
whence, in particular, $ |\nabla {\bf v}|^{p-2} \nabla {\bf v}  \in L^{p'}(\o , \rNn)$. Therefore, by a density argument, equation \eqref{main101bis} also holds for every $\bm{\varphi} \in W^{1,p}_0(\o, \rN)$. This means that   ${\bf v}$ is a weak solution to the Dirichlet problem \eqref{localeq}$+$\eqref{dircond}. Its 
uniqueness  ensures that  ${\bf v}=\bu$.
\color{black}
Furthermore, owing to \eqref{main27aus} and \eqref{main32}, 
\begin{equation}\label{main35}
\||\nabla \bu|^{p-2} \nabla \bu\|_{W^{1,2}(\o, \rNn)} \leq C \|{\bf f}\|_{L^2(\o, \rN)}
\end{equation}
%
%
for some constant $C=C(n, N, p, L_\o, d_\o, \mathcal K_\o)$.

\smallskip
\par\noindent
\emph{Step 3}. We conclude by removing the remaining additional assumption \eqref{fsmooth}. Let ${\bf f} \in L^2(\Omega, \rN)$. By the definition of approximable solution, there exists a  sequence
$\{{\bf f}_k\} \subset C^\infty_0(\Omega, \rN )$, with  ${\bf f}_k \to {\bf f}$ in $L^2(\o, \rN)$,  such that the sequence of weak solutions 
 $\{\bu _k\} \subset W^{1,p}_0(\o, \rN)$ to problems \eqref{eqdirichletk}, satisfies 
  $\bu _k \to \bu$ and $\nabla \bu_k \to \nabla \bu$ a.e. in $\o$.
%
%
%
%
By inequality \eqref{main35} of the previous step, applied with $\bm u$ and $\bm f$ replaced by $\bm u _k$ and $\bm f_k$, we have that $|\nabla \bu_k|^{p-2} \nabla \bu_k \in W^{1,2}(\Omega, \rNn)$, and there exist  constants $C_1$ and $C_2$, depending on $N$, $p$ and $\o$, such that
\begin{align}\label{main40}
\||\nabla \bu_k|^{p-2} \nabla \bu_k\|_{W^{1,2}(\o, \rNn)} \leq C_1 \|{\bf f}_k\|_{L^2(\o, \rN)} \leq C_2 \|{\bf f}\|_{L^2(\o, \rN)}\,.
\end{align}
Hence, the sequence $\{|\nabla \bu_k|^{p-2}\nabla \bu_k\}$ is uniformly bounded in $ W^{1,2}(\Omega, \rNn)$, and there exists a subsequence, still indexed by $k$, and a function ${\bf U} \in W^{1,2}(\Omega, \rNn)$ fulfilling
\begin{equation}\label{main38}
|\nabla \bu_{k}|^{p-2}\nabla \bu_{k} \to {\bf U}\quad \hbox{in $L^2(\Omega, \rNn)$,} \quad  |\nabla \bu_{k}|^{p-2} \nabla \bu_{k} \rightharpoonup {\bf U}\quad \hbox{in $W^{1,2}(\Omega, \rNn)$}.
\end{equation}
Since  $\nabla \bu_k \to \nabla \bu$ a.e. in $\o$, we thus infer that $|\nabla \bu|^{p-2} \nabla \bu ={\bf U}  \in W^{1,2}(\Omega, \rNn)$,
and the second inequality in \eqref{seconddir2} follows via \eqref{main40} and     \eqref{main38}.  The first inequality  in \eqref{seconddir2} holds trivially.
%
%
The statement concerning the solution to  the Dirichlet problem \eqref{localeq}$+$\eqref{dircond} is thus fully proved. 

\smallskip
\par\noindent
The outline of the proof  for  the solution to the Neumann problem \eqref{localeq}$+$\eqref{neucond} is the same as that for the Dirichlet problem. We point out hereafter just the variants required in the various steps.
\\ \emph{Step 1}. The solution $\bm u_\ep$ to the Dirichlet  problem \eqref{eqdirichletep} must be replaced, of course, by a solution $\bm u_\ep$ to the approximating Neumann problem
\begin{equation}\label{eqneumannep}
\begin{cases}
-{\rm {\bf  div}} (b_\varepsilon  (|\nabla \bu_\varepsilon|)^{p-2} \nabla \bu_\varepsilon) = {\bf f} & {\rm in}\,\,\, \o \\
 \displaystyle  \frac{\partial \bu_\varepsilon}{\partial \bm \nu} =0  &
{\rm on}\,\,\,
\partial \o \,.
\end{cases}
\end{equation}
Such a 
 solution   is only unique up to additive constant vectors. 
%
Estimate \eqref{gradbound} is a consequence of \cite[Theorem 2.4]{cmARMA}. As shown in the proof of that theorem, $\bm u_\ep \in W^{1,\infty}(\o, \rN) \cap W^{2,2}(\o, \rN)$, and there exists 
a
sequence $\{{\bf u}_k\} \subset C^\infty (\Omega, \rN)\cap
C^2(\overline \Omega, \rN)$ enjoying the following properties: $$ \frac{\partial \bu_k}{\partial \bm \nu} = 0 \quad \hbox{ on $\partial
\Omega$\,,}$$
\begin{equation}\label{convuk}
{\bf u}_k \to {\bf u}_\varepsilon \quad \hbox{in $W^{2,2}(\Omega, \rN)$,}
\quad \nabla {\bf  u}_k \to \nabla {\bf u}_\varepsilon 
\quad \hbox{a.e. in
$\Omega $},
\end{equation}
\begin{equation}\label{gradboundk}
\|\nabla \bu_k\|_{L^\infty (\o, \rNn)}  \leq C \|\nabla \bu_\varepsilon\|_{L^\infty (\o, \rNn)} 
\end{equation}
\color{black}
for some constant $C$ independent of $k$, 
%
%
 $|\nabla |\nabla {\bf u}_k|| \leq |\nabla ^2
\bu_k|$ a.e. in $\Omega$,
and 
\begin{align}\label{sys1}
- {\rm {\bf div}} (b_\varepsilon(|\nabla {\bf u_k}|)^{p-2}\nabla {\bf u_k} ) \to {\bf
f}
\quad \hbox{in $L^2(\Omega , \rN)$},
\end{align}
as $k \to \infty$. Furthermore, equation \eqref{napoli} now holds owing to \eqref{neumanngrad}, with a constant $C=C(\o, N, p)$ depending on $\o$ via an upper bound for $d_\o$ and $L_\o$. With these variants in place,   inequality \eqref{main27} follows via the same argument as in the case of the Dirichlet problem.
%
%
%

\smallskip
\par\noindent
 \emph{Step 2}. 
The Dirichlet problem \eqref{eqm} has to be replaced with the Neumann problem
\begin{equation}\label{eqmneu}
\begin{cases}
- {\rm {\bf div}} (|\nabla \bu_m|^{p-2} \nabla \bu_m ) = {\bf f}  & {\rm in}\,\,\, \o _m \\
 \displaystyle \frac{\partial \bu_m}{\partial \bm \nu} =0  &
{\rm on}\,\,\,
\partial \o _m \,,
\end{cases}
\end{equation}
and the corresponding sequence of solutions $\{\bu_m\}$ has  to be normalized by a suitable sequence of additive constant vectors 
$\bm \xi _m \in \rN$ in such a way that inequality \eqref{energym} still holds. For instance, one can choose $\bm \xi _m = - \int _\o \bu _m\, dx$. After this normalization, the sequence $\{\bu_m\}$ admits a subsequence, still denoted by $\{\bu_m\}$, converging to a function $\bm v$ with the same properties as in the case of the Dirichlet problem.
\\ Next, any  test function $\bm{\varphi} \in W^{1, \infty}(\o, \rN)$ can be exended to a function in $W^{1, \infty}(\mathbb R^n, \rN)$, still denoted by $\bm{\varphi}$. 
The use of this test function in the weak formulation of problem \eqref{eqmneu} yields
\begin{equation}\label{main101neumann}
\int _{\Omega_m} |\nabla \bu_m|^{p-2} \nabla \bu_m \cdot \nabla \bm{\varphi} \, dx = \int_{\Omega _m} {\bf f} \cdot \bm{\varphi} \, dx\,.
\end{equation}
Passage to the limit as $m \to \infty$ on the left-hand side of equation \eqref{main101neumann}, to obtain 
\begin{equation}\label{main101limit}
\int _{\Omega} |\nabla \bm v|^{p-2} \nabla \bm v \cdot \nabla \bm{\varphi} \, dx = \int_{\Omega} {\bf f} \cdot \bm{\varphi} \, dx\,,
\end{equation}
can be justified as follows.  The left-hand side of equation \eqref{main101neumann} can be split as 
\begin{equation}\label{main101'}
\int _{\o _m} |\nabla \bu_m|^{p-2}\nabla \bu_m \cdot \nabla  \bm{\varphi} \, dx = 
\int _{\o } |\nabla \bu_m|^{p-2}\nabla \bu_m \cdot \nabla \bm{\varphi}  \, dx + \int _{\o _m \setminus \o} |\nabla \bu_m|^{p-2}\nabla \bu_m \cdot \nabla \bm{\varphi}  \, dx\,.
\end{equation}
The first integral on the right-hand side of \eqref{main101'} converges to the left-hand side of \eqref{main101limit} as $m \to \infty$, 
%
%
%
owing to \eqref{main32} and \eqref{main100}. The second integral tends to $0$, by \eqref{main27aus} and the fact that $|\o _m \setminus \o| \to 0$. 
\\ Since $\o$ is a bounded Lipschitz domain, the density of $W^{1, \infty}(\o, \rN)$ in $W^{1, p}(\o, \rN)$ entails that equation \eqref{main101limit} continues to hold for any test function $\bm{\varphi}$ in the latter space, and hence $\bm v$ is the (unique, up to additive constant vectors) weak solution to the Neumann problem \eqref{localeq}$+$\eqref{neucond}.

\smallskip
\par\noindent
\emph{Step 3}. This step is completely analogous, provided that the sequences $\{\bm{f}_k\}$ and $\{\bm{u}_k\}$ are taken as in the definition of approximable solution $\bm u$ to the Neumann problem  \eqref{localeq}$+$\eqref{neucond} .
%
%
\qed

\medskip
\par\noindent
{\bf Proof of Corollary  \ref{secondmarc}}.  Lemmas \ref{tracecap} and \ref{traceineq} ensure that
\begin {equation}\label{gen3} \mathcal K_\o (r) \leq C \sup_{x\in \partial \o}\|\mathcal B\|_{X(\partial \o \cap B_r(x))} \qquad \hbox{for $r \in (0, r_0)$,}
\end{equation}
for suitable constants $r_0$ and $C$ depending on $n$, $L_\o$ and $d_\o$.
The conclusion then follows from Theorem \ref{seconddir},  via inequality \eqref{gen3}. \qed

\medskip
\par\noindent
{\bf Proof of Theorem \ref{secondconvex}}.  The proof is analogous to that of Theorem \ref{seconddir}. The only difference is that, in Step 2, one has  to  choose a sequence $\{\Omega _m\}$ of bounded convex open sets,  with $\partial \o _m \in C^\infty$, approximating $\o$ from outside with respect to the Hausdorff distance. Conditions \eqref{feb35} are authomatically fulfilled in this case. On the other hand, condition \eqref{appcap0} is irrelevant in the present situation, thanks to 
 the fact that the constant $C$ in inequality \eqref{fund} is independent of the function $\mathcal K_\o$ in the case of convex domains $\o$.
\qed

\section*{Appendix. 
}

Here, we present an outline of the proof of Lemma \ref{approxcap}. 


\par\noindent
Let
 $\{\sigma _m\}$ be a sequence of radially symmetric mollifiers in $\mathbb R^{n-1}$,
namely, $\sigma _m \in C^\infty_0(\mathbb R^{n-1})$, ${\rm supp}\, \sigma _m \subset B_{1/m}^{n-1}(0)$, $\sigma _m \geq 0$ and $\int _{\mathbb R^{n-1}}\sigma _m\, dx'=1$ for $m \in \mathbb N$. Here, the index $n-1$ attached to the notation of a ball denotes a ball in $\mathbb R^{n-1}$. 
Given $g\in L^1_{\rm loc}(\mathbb R^{n-1}, \mathbb R^d)$, with $d \in \mathbb N$, we denote by $M_m(g)$ the convolution of $g$ with $\sigma_m$, namely the function  $M_m(g) : \rn \to \mathbb R^d$ defined as 
$$M_m(g) (x')= \int _{\mathbb R^{n-1}} g(y') \sigma _m(x'-y')\, dy' \quad \hbox{for $x' \in \mathbb R^{n-1}$.}$$
Moreover, given a function $h \in L^1_{\rm loc}(\rn,  \mathbb R^d)$,  define $\widetilde M_m(h) : \rn \to \mathbb R^d$ as the convolution of $h$ with $\sigma _m$ with respect to the first $n-1$ variables, i.e.
$$\widetilde M_m(h)(x', x_n) = \int _{\mathbb R^{n-1}} h(y',x_n) \sigma _m(x'-y')\, dy' \quad \hbox{for $(x',x_n) \in \rn$.}$$
Note that, if ${\rm supp} \, h \subset \subset B_r$ for some ball $B_r\subset \rn$, then ${\rm supp} \widetilde M_m(h) \subset \subset B_r$ as well, if $m$ is sufficiently large.
\\
Assume, for the time being, that  $\o\subset \rn$ globally agrees with the subgraph of a Lipschitz continuous function $\phi : \mathbb R^{n-1} \to \mathbb R$,    let $\varrho$ be a nonnegative function in $L_{\rm loc}^1(\partial \o)$, and let $r_0\in (0,1)$. Define, for $x \in \partial \o$ and $r \in (0,r_0)$,
\begin{equation}\label{main6}
R_{B_r(x)}(\varrho, \o)  = \sup_{v \in C^{0,1}_0(B_r(x))} \frac{\int _{\partial \o \cap B_r(x)} v^2\, \varrho d\hh}{\int _{B_r(x)}|\nabla v|^2 \, dx}\,,
%
\end{equation}
and 
\begin{equation}\label{tracenmu}
Q_{B_r(x)}(\varrho, \o) = \sup_{E\subset \partial\o \cap B_r(x)}
 \frac{\int _{E}\varrho \, d\hh}{{\rm cap}_{B_1(x)}(E)}\,.
\end{equation}
Then, by Lemma \ref{tracecap} and Remark \ref{pallepiene},
\begin{equation}\label{appcap7}
R_{B_r(x)}(\varrho, \o) \approx Q_{B_r(x)}(\varrho, \o)\,,
\end{equation}
up to positive multilicative constants depending only on $n$, $r_0$ and on an upper bound for the Lipschitz constant $L$ of $\phi$.
\\
Assume, in addition, that   $\phi \in W^{2,1}(\mathbb R^{n-1})$, and let $\mathcal B_\phi$ denote the (weak) second fundamental form of the graph of $\phi$.
Then,
\begin{equation}\label{appcap5}
\mathcal B_\phi  = \frac {\nabla ^2\phi}{\sqrt{1+ |\nabla \phi|^2}}\,.
\end{equation}
Hence,
\begin{equation}\label{secondform}
\frac{|\nabla ^2 \phi (x')|}{\sqrt{1+ L^2}}  \leq |\mathcal B_\phi (x')| 
\leq |\nabla ^2 \phi (x')|  \quad \hbox{for a.e. $x' \in \mathbb R^{n-1}$,}
\end{equation}
and
\begin{align}\label{appcap1} 
 |\mathcal B_{M_m(\phi)}(x')| \leq  |\nabla ^2 M_m(\phi)(x')|= |M_m(\nabla ^2 \phi)(x')| \leq M_m(|\nabla ^2 \phi|)(x')\quad  \hbox{for a.e. $x' \in \mathbb R^{n-1}$.}
\end{align}
 Given $v\in C^{0,1}_0(B_r(0))$ for some $r\in (0,r_0)$,
define the function $w_m: \rn \to \mathbb R$ as
$$w_m= \sqrt{\widetilde M_m(v^2)}\,.$$
As noticed above,   
\begin{equation}\label{suppwk}
{\rm supp}\, w_m \subset \subset B_r(0)
\end{equation}
 if $m$ is sufficiently large, owing to our choice of $v$.
Moreover, on denoting by $\nabla _{x'}$ the gradient operator with respect to the sole variables $x'\in \mathbb R^{n-1}$, one has that
\begin{equation}\label{appcap3}
\big|\nabla _{x'}\big(\widetilde M_m(v^2)\big)\big| = \big| \widetilde M_m\big(\nabla _{x'} (v^2)\big)\big|= 2  \big|( \widetilde M_m\big(v \nabla _{x'}v)\big)\big|\leq 2 \sqrt{\widetilde M_m(v^2)}\, \sqrt{\widetilde M_m\big(|\nabla _{x'}v|^2\big)}\,.
\end{equation}
Also, 
\begin{equation}\label{appcap4}
\big|\big(\widetilde M_m(v^2)\big)_{x_n}\big| = \big|\widetilde M_m\big((v^2)_{x_n}\big)\big|= 2  \big|(\widetilde M_m(v v_{x_n}))\big|\leq 2 \sqrt{\widetilde M_m(v^2)}\, \sqrt{\widetilde M_m\big((v_{x_n})^2\big)}.
\end{equation}
Thus,  $w_m$ is Lipschitz continuous, \color{black} and 
\begin{equation}\label{appcap5bis}
|\nabla w_m| \leq C \frac{\big|\nabla _{x'}\big(\widetilde M_m(v^2)\big)\big| + \big|\big(\widetilde M_m(v^2)\big)_{x_n}\big|}{\sqrt{\widetilde M_m(v^2)}} \leq  C' \sqrt{\widetilde M_m\big(|\nabla _{x'}v|^2\big)   + \widetilde M_m\big((v_{x_n})^2\big)} \quad \hbox{a.e. in $\mathbb R^{n-1}$}
\end{equation}
for some absolute constants $C$ and $C'$. Hence, via an application of Fubini's theorem, we deduce that
\begin{align}\label{appcap9}
\int _{\rn} |\nabla w_m|^2\, dx& \leq C' \int _{\rn} \widetilde M_m\big(|\nabla _{x'}v|^2\big)   + \widetilde M_m\big((v_{x_n})^2\big)\, dx \\ \nonumber & = C' \int _{\rn}  |\nabla _{x'}v|^2   +  (v_{x_n})^2\, dx
 =C' \int _{\rn}  |\nabla  v|^2  \, dx\,.
\end{align}
Owing to equations \eqref{secondform} and \eqref{appcap1},
\begin{align}\label{appcap2}
\int _{\mathbb R^{n-1}} v^2(x', 0) \, |\mathcal B_{M_m(\phi)}(x')|\, d x'  & \leq 
\int _{\mathbb R^{n-1}} v^2(x', 0) \, M_m\big(|\nabla ^2 \phi (x')|\big)\, d x' \\ \nonumber & = 
\int _{\mathbb R^{n-1}} \widetilde M_m(v^2(x', 0)) \, |\nabla ^2 \phi (x')|\, d x' 
\\ \nonumber & \leq\sqrt{1+ L^2}   \int _{\mathbb R^{n-1}} \widetilde M_m(v^2(x', 0)) |\mathcal B_{\phi} (x')|\, d x
\\ \nonumber &
= \sqrt{1+ L^2} \int _{\mathbb R^{n-1}} w_m(x',0)^2 \,  |\mathcal B_{\phi} (x')|\, d x' \,.
\end{align}
By definition \eqref{main6} and equation \eqref{suppwk},
\begin{equation}\label{appcap30}
\frac{\int _{\mathbb R^{n-1}} w_m(x',0)^2 \,  |\mathcal B_{\phi} (x')|\, d x'}{\int _{\rn} |\nabla w_m|^2\, dx} \leq R_{B_r(0)}(|\mathcal B_{\phi}|, \mathbb R^n _-)\,,
\end{equation}
provided that $m$ is sufficiently large, where we have set $\mathbb R^n _- = \{(x',x_n): x_n<0\}$.
Inequalities \eqref{appcap9}--\eqref{appcap30} imply that
%
%
%
\begin{align}\label{appcap8}
R_{B_r(0)}(|\mathcal B_{M_m(\phi)}|, \mathbb R^n _-)    = \sup_{v \in C^{0,1}_0(B_r(0))} \frac{\int _{\mathbb R^{n-1}} v^2(x', 0) \, |\mathcal B_{M_m(\phi)}(x')|\, d x'}{\int _{\rn}|\nabla v|^2 \, dx}  \leq C R_{B_r(0)}(|\mathcal B_\phi|, \mathbb R^n _-)\,,
\end{align}
for some constant $C=C(n, L, r_0)$.
\\ Now, denote by $G_\phi$ the graph of $\phi$ and by $S_\phi$ the subgraph of $\phi$, and define $G_{M_m(\phi)}$ and $S_{M_m(\phi)}$ analogously. Set $x^0= (0, \phi (0))\in \rn$ and $x_m^0=(0, M_m(\phi)(0))\in \rn$.
Given $v \in C^{0,1}_0(B_r(x_m^0))$, we have that
\begin{align}\label{appcap10}
\int _{G_{M_m(\phi)}} v^2|\mathcal B_{M_m(\phi)}| d\hh & = \int _{\mathbb R^{n-1}} v^2(x', M_m(\phi)(x')) \, |\mathcal B_{M_m(\phi)}(x')|\, \sqrt{1+ |\nabla M_m(\phi)|^2} \,d x' \\
\nonumber & \leq 
 \sqrt{1+ L^2} \int _{\mathbb R^{n-1}} v^2(x', M_m(\phi)(x')) \, |\mathcal B_{M_m(\phi)}(x')|\,  \,d x'\,.
\end{align}
Define the function $w_m : \rn \to \mathbb R$ as
\begin{equation}\label{appcap34}
w_m(x',x_n) = v(x', x_n+ M_m(\phi)(x')) \quad \hbox{for $(x',x_n) \in \rn$.}
\end{equation}
Then, ${\rm supp} \,w_m \subset \subset B_r(0)$ if $m$ is sufficiently large. Furthermore, 
\begin{align}\label{appcap11}
\int _{\mathbb R^{n-1}} v^2(x', M_m(\phi)(x')) \, |\mathcal B_{M_m(\phi)}(x')|\,d x' = \int _{\mathbb R^{n-1}} w^2_m(x', 0) \, |\mathcal B_{M_m(\phi)}(x')|\,d x'\,.
\end{align}
By \eqref{appcap8} (with $v$ replaced by $w_m$) and \eqref{appcap9},
\begin{equation}\label{appcap12}
\int _{\mathbb R^{n-1}} w^2_m(x', 0) \,|\mathcal B_{M_m(\phi)}(x')|\,d x' \leq C R_{B_r(0)}(|\mathcal B_\phi|,  \mathbb R^n_-)
\int _{\mathbb R^{n}} |\nabla w_m|^2\, dx\,.
\end{equation}
On the other hand, since 
$$|\nabla  M_m(\phi)| = |M_m(\nabla  \phi)| \leq L \quad \hbox{a.e. in $\mathbb R^{n-1}$\,,}$$
for $m \in \mathbb N$, there exists a constant $C=C(n,L)$ such that
\begin{equation}\label{appcap13}
|\nabla w_m(x', x_n)| \leq C |\nabla v(x', x_n+ M_m(\phi)(x'))| \quad \hbox{for a.e. $(x',x_n) \in \rn$.}
\end{equation}
Thereby, 
\begin{equation}\label{appcap14}
\int _{\mathbb R^{n}} |\nabla w_m|^2\, dx \leq C^2 \int _{\mathbb R^{n}}  |\nabla v(x', x_n+ M_m(\phi)(x'))|^2\, dx = C^2 \int _{\mathbb R^{n}}  |\nabla v|^2\, dx \,.
\end{equation}
Combining equations \eqref{appcap10}, \eqref{appcap11}, \eqref{appcap12} and \eqref{appcap14} yields
\begin{equation}\label{appcap15}
\int _{G_{M_m(\phi)}} v^2|\mathcal B_{M_m(\phi)}| d\hh \leq C R_{B_r(0)}(|\mathcal B_\phi|, \mathbb R^n_-) \int _{\mathbb R^{n}}  |\nabla v|^2\, dx 
\end{equation}
for some constant $C=C(n,L, r_0)$ and every function $v \in C^{0,1}_0(B_r(x_m^0))$. Therefore, 
\begin{equation}\label{gen1}
R_{B_r(x_m^0)}(|\mathcal B_{M_m(\phi)}|, S_{M_m(\phi)}) \leq C R_{B_r(0)}(|\mathcal B_\phi|, \mathbb R^n_-) 
\end{equation}
for some constant $C=C(n,L, r_0)$.
On the other hand, an analogous change of variable as in \eqref{appcap34}, with $M_m(\phi)$ replaced by $- \phi$, tells us that
\begin{equation}\label{appcap35rev}
R_{B_r(0)}(|\mathcal B_\phi|, \mathbb R^N_-) 
\leq C 
R_{B_r(x^0)}(|\mathcal B_{\phi}|, S_{\phi})
%
\end{equation}
for some constant $C=C(n,L, r_0)$.
Inequalities \eqref{gen1} and \eqref{appcap35rev} yield
\begin{equation}\label{appcap16}
R_{B_r(x_m^0)}(|\mathcal B_{M_m(\phi)}|, S_{M_m(\phi)}) \leq C  R_{B_r(x^0)}(|\mathcal B_{\phi}|, S_{\phi})
\end{equation}
for some constant $C=C(n,L, r_0)$, for every sufficiently large $m \in \mathbb N$. Of course, inequality \eqref{appcap16} continues to hold if 
$x^0=(x', \phi (x'))$ and $x_m^0= (x', M_m(\phi)(x'))$   for any $x'\in \mathbb R^{n-1}$.
\\
Assume now that $\o$ is a bounded Lipschitz domain in $\rn$. Let 
  $\mathcal U= B_\rho^{n-1}\times (-\delta , \delta)$ for some $\rho>0$ and $\delta >0$, and let $\phi \in W^{2,1}(\mathbb R^{n-1})$ is a Lipschitz continuous function, with ${\rm supp}\, \phi \subset  B_\rho^{n-1}$, such that (on translating and rotating $\o$, if necessary) 
$\Omega \cap \mathcal U   = \{(x',x_n): x' \in  B_\rho^{n-1}, -\delta < x_n< \phi(x')\}$ and $ \partial \Omega \cap \mathcal U   = \{(x',x_n): x' \in  B_\rho^{n-1}, x_n= \phi(x')\}$.
Let  $M_m(\phi)$ be defined as above for $m \in \mathbb N$. Denote by $\Omega _m$ the open subset which agrees with $\o$ outside $\mathcal U$, and satisfying 
$
\Omega_m \cap \mathcal U  = \{(x',x_n): x' \in  B_\rho^{n-1}, -\delta < x_n< M_m(\phi)(x')\}$ and  $\partial \Omega_m \cap \mathcal U  = \{(x',x_n): x' \in  B_\rho^{n-1}, x_n= M_m(\phi)(x')\}$.
From inequality \eqref{appcap16} one can deduce that
\begin{align}\label{appcap22}
R_{B_r(x_m^0)}(|\mathcal B_{M_m(\phi)}|, \o_m)  \leq C R_{B_{r}(x^0)}(|\mathcal B_\phi|,  \o),
\end{align}
for every ball $B_r(x_m^0)\subset \mathcal U$.
Since $M_m(\phi) \to \phi$ as $m \to \infty$, inequality \eqref{appcap22} can be used, via \eqref{appcap7}, to prove inequality
\eqref{appcap0} when $\{\o_m\}$ is a   sequence  of sets obtained from a local smoothening of $\partial \o$. 
\\ The construction of a sequence $\{\o_m\}$ of sets as in the statement can be accomplished on combining this local construction with a partition of unit argument. For brevity, we limit ourselves to  sketching this argument. To begin with, note that the sequence of functions $\{M_m(\phi)\}$ above can be modified, if necessary, by an additive constant in such a way the modified sequence satisfies $\phi \leq M_m(\phi)$ in $B_r^{n-1}$, still converges to $\phi$ and fulfills inequality \eqref{appcap22}.
%
Next, since $\o$ is a bounded Lipschitz domain with $\partial \o \in W^{2,1}$, there exists a  finite family of open sets $\{\mathcal U_i\}$, and a corresponding family  $\{\mathcal \phi _i\}$ of Lipschitz continuous  functions with $\phi _i \in W^{2,1}(\mathcal U_i)$, such that 
\begin{equation}\label{appcap40}
\partial \Omega =  \bigcup _i\mathcal V_i\,,
\end{equation}
where $\mathcal V_i$ are open subsets in $\partial \Omega$ (with respect to the topology induced by $\rn$), such that
\begin{equation}\label{appcap41}
\mathcal V_i = \Phi _i (\mathcal W_i)\,,
\end{equation}
$\mathcal W_i$ is an open subset of $\mathbb R^{n-1}$, and
\begin{equation}\label{appcap42}
\Phi _i : \mathcal W_i \to \mathcal V_i
\end{equation}
is a homeomerphism obtained by composing the function 
$$ \mathcal U_i \ni x' \mapsto (x', \phi _i(x')) \in \rn$$
with an isometry. Let $M_m(\phi _i)$ be the function defined as above, with $\phi$ replaced by $\phi_i$, and let $\Phi _i^m: \mathcal W_i \to \rn$ be the function obtained on replacing $\phi_i$ with $M_m(\phi _i)$ in  definition \eqref{appcap42} of $\Phi _i$. 
Let $\{\xi _i\}$ be a partition of unity associated with the coverning $\{\mathcal V_i\}$ of $\partial \Omega$. Define, for $m \in \mathbb R$, the function $\Psi _m : \partial \Omega \to \rn$ as
$$\Psi _m = \sum _i \big(\Phi_i^m \circ \Phi_i^{-1}\big)\xi _i\,.$$
Then $\o_m$ can be defined as the open set in $\rn$ such that  
$\partial \o _m = \Psi_m(\partial \Omega)$.
The fact that the sequence $\{\o_m\}$ fulfills property \eqref{appcap0} can be deduced as a  consequence of the local argument described above. The remaining properties described in the statement are easier to verify. 

\bigskip
\par\noindent
{
\color{black}
{\bf Acknowledgment}. This research was partly carried out during a visit of the first-named author at the Institut Mittag-Leffler in August 2017. He wishes to thank the Director and the staff of the Institut for their support and hospitality.}

\section*{Compliance with Ethical Standards}\label{conflicts}

\smallskip
\par\noindent 
{\bf Funding}. This research was partly funded by:   
\\ (i) Research Project of the
Italian Ministry of University and Research (MIUR) Prin 2012 \lq\lq Elliptic and
parabolic partial differential equations: geometric aspects, related
inequalities, and applications"  (grant number 2012TC7588);    
\\ (ii) GNAMPA   of the Italian INdAM - National Institute of High Mathematics (grant number not available);   
\\  (iii) Ministry of Education and Science of the Russian Federation  (grant number 02.a03.21.0008).

\smallskip
\par\noindent
{\bf Conflict of Interest}. The authors declare that they have no conflict of interest.


\begin{thebibliography}{99}


 




%

\bibitem[ACMM]{ACMM} A.Alvino, A.Cianchi, V.G.Maz'ya \&
A.Mercaldo, Well-posed elliptic Neumann problems involving
irregular data and domains,  \emph{Ann. Inst. H. Poincar\'e Anal. Non Lin\'eaire} {\bf 27}  (2010), 1017--1054.


\bibitem[AKM]{AKM} B.Avelin, T.Kuusi \&
G.Mingione, Nonlinear Calder\'on-Zygmund theory in a limiting case,  \emph{Arch. Rat. Mech. Anal.} {\bf 227} (2018),  663--714.



%
%




\bibitem[BBGGPV]{BBGGPV} P.B\'enilan, L.Boccardo, T.Gallou\"et, R.Gariepy, M.Pierre \& J.L.Vazquez,
An $L^1$-theory of existence and uniqueness of solutions of
nonlinear elliptic equations, \emph{Ann. Sc. Norm. Sup. Pisa} {\bf
22} (1995), 241--273.

\bibitem[BF]{BF} A.Bensoussan \& J.Frehse,  ``Regularity results for nonlinear elliptic systems and applications", Springer-Verlag, Berlin, 2002.

\bibitem[BeCr]{BC} H.Beir\~{a}o da Veiga \& F.Crispo, On the global $W^{2,q}$ regularity for
nonlinear $N$-systems of the $p$-Laplacian type in $n$ space
variables, \emph{Nonlinear Anal.} {\bf 75} (2012), 4346--4354.

\bibitem[Be]{Bernstein} S.Bernstein, Sur la nature analytique des solutions des \'equations aux d\'eriv\'ees partielles du second ordre, 
{\em Math. Ann.}  {\bf 59} (1904), 20--76


\bibitem[BoGa]{BG} L.Boccardo \& T.Gallou\"et, Nonlinear elliptic and parabolic equations involving measure data,
  \emph{J. Funct. Anal.} {\bf 87} (1989),
 149--169.

%


\bibitem[BCDKS]{BCDKS} D.Breit, A.Cianchi, L.Diening, T.Kuusi \& S.Schwarzacher, Pointwise Calder\'on-Zygmund gradient estimates for the $p$-Laplace system,
\emph{J. Math. Pures. Appl.}, to appear.


%

%

 \bibitem[Bu]{Burenkov} V.Burenkov,  ``Sobolev spaces on domains", Teubner, Stuttgart, 2002.


 \bibitem[Ce1]{Cel} A.Cellina, The regularity of solutions to some variational problems, including the $p$-Laplace equation, for $2 \leq p <3$,  \emph{ESAIM COCV} {\bf 23} (2017), 1543--1553.



\bibitem[Ce2]{Cel1} A.Cellina, The regularity of solutions to some variational problems, including the $p$-Laplace equation, for $3 \leq p <4$, \emph{preprint}.




\bibitem[ChDi]{CDiB} Y.Z.Chen \& E.Di Benedetto, Boundary estimates
for solutions of nonlinear degenerate parabolic systems,  \emph{J.
Reine Angew. Math.} {\bf 395} (1989), 102--131

\bibitem[CiMa1]{cmARMA} A.Cianchi \& V.Maz'ya, Global boundedness of the gradient for a class of nonlinear elliptic systems,
  \emph{Arch. Ration. Mech. Anal.} {\bf 212} (2014), 129--177.

\bibitem[CiMa2]{cmsecond} A.Cianchi \& V.Maz'ya,  Second-order two-sided estimates in nonlinear elliptic problems,  \emph{Arch. Ration. Mech. Anal.}, {\bf 229} 
  (2018),  569--599.
%

\bibitem[CGM]{CGM} F.Crispo, C.R.Grisanti \& P.Maremonti, On the high regularity of solutions to the $p$-Laplacian boundary value problem in exterior domains, \emph{Ann.
Mat. Pura Appl.}  {\bf 195}  (2016), 821--834.


\bibitem[DaA]{DallA} A.Dall'Aglio, Approximated
solutions of equations with $L\sp 1$ data. Application to the
$H$-convergence of quasi-linear parabolic equations, \emph{Ann.
Mat. Pura Appl.}  {\bf 170}  (1996), 207--240.



\bibitem[DaSc]{Dasc}
L.Damascelli \& B.Sciunzi, Regularity, monotonicity and symmetry of positive solutions of $m$-Laplace equations,   \emph{J. Diff. Equat.}  {\bf 206}  (2004), 483--515.

%


%
%
%
%


\bibitem[DHM1]{DHM1} G.Dolzmann, N.Hungerb\"uhler \& S.M\"uller, The $p$-harmonic system with measure-valued right hand side,  \emph{Ann. Inst. H. Poincar\'e Anal. Non Lin\'eaire} {\bf 14} (1997),   353--364. 
%

 \bibitem[DHM2]{DHM3} G.Dolzmann, N.Hungerb\"uhler \& S.M\"uller,  Uniqueness and maximal regularity for nonlinear elliptic systems of $n$-Laplace type with measure valued right hand side,  \emph{J. Reine Angew. Math.}  {\bf 520} (2000), 1--35. 

%



\bibitem[DuMi1]{DM} F.Duzaar \& G.Mingione, Gradient estimates via non-linear potentials, \emph{Amer. J. Math.}  {\bf 133} (2011), 1093--1149.





\bibitem[Gi]{giaquinta} M.Giaquinta,  ``Multiple integrals in the calculus of variations and nonlinear elliptic systems", Annals of Mathematical Studies,
 Princeton University Press, Princeton, NJ, 1983.





\bibitem[Gr]{Grisvard} P.Grisvard,  ``Elliptic problems in nonsmooth domains", Pitman, Boston, MA, 1985.





%



%








\bibitem[KrMa]{KrolM} I.N.Krol' \& V.G. Maz'ya, On the absence
of continuity and H\"older continuity of solutions of quasilinear
elliptic equations near a nonregular boundary,  \emph{Trudy
Moskov. Mat. Os\v{s}\v{c}.} {\bf 26} (1972) (Russian); English
translation: \emph{Trans. Moscow Math. Soc.} {\bf 26} (1972),
73--93.









\bibitem[LiMu]{LM}  P.-L.Lions \& F.Murat, Sur les
solutions renormalis\'ees d'\'equations elliptiques non
lin\'eaires, manuscript.

\bibitem[Lo]{Lou} H.Lou, 
On Singular Sets of Local Solutions to $p$-Laplace Equations,  \emph{Chin. Ann. Math.} {\bf 29B} (2008),   521--530.




\bibitem[MaMi]{MarcusMizel} M.Marcus \& V.J.Mizel,  Absolute continuity of tracks and mappings of Sobolev
spaces, \emph{Arch. Ration. Mech. Anal.} {\bf 45}
 (1972),  294--320.






\bibitem[Maz1]{Mazya62} V.Maz'ya, The negative spectrum of the higher-dimensional Schr\"odinger operator,  {\em Dokl. Akad. Nauk SSSR} {\bf 144} (1962), 721--722 (Russian). English translation: \emph{Sov. Math. Dokl.}  {\bf 3} (1962). 

\bibitem[Maz2]{Mazya64} V.Maz’ya, On the theory of the higher-dimensional Schr\"odinger operator,  {\em Izv. Akad. Nauk SSSR} Ser. Mat. {\bf 28} (1964), 1145--1172 (Russian).




 \bibitem[Maz3]{Ma67} V.G.Maz'ya, Solvability in $W^2_2$ of the Dirichlet problem in a region with a smooth irregular boundary,
{\em Vestnik Leningrad. Univ.} {\bf 22} (1967), 87--95 (Russian).

 \bibitem[Maz4]{Ma73} V.G.Maz'ya, The coercivity of the Dirichlet problem in a domain with irregular boundary,
{\em Izv. Vys\v s. U\v cebn. Zaved. Matematika} {\bf 4} (1973), 64--76 (Russian).

%


\bibitem[Maz5]{Ma69} V.G.Maz'ya, On weak solutions of the Dirichlet and Neumann problems,
{\em Trusdy Moskov. Mat. Ob\v s\v c.} {\bf 20} (1969), 137--172
(Russian); English translation: \emph{Trans. Moscow Math. Soc.}
{\bf 20} (1969), 135--172.
\bibitem[Maz6]{Mabook} V.G.Maz'ya,  ``Sobolev spaces with applications to elliptic partial differential equations", Springer-Verlag, Heidelberg, 2011.
\bibitem[MazSh]{MaSh} V.G. Maz'ya \& T.O.Shaposhnikova, \lq\lq Theory of Sobolev multipliers. With applications to differential and integral operators." Springer-Verlag, Berlin, 2009.


\bibitem[Mi]{Mi1} G.Mingione, Gradient estimates below the duality
exponent, \emph {Math. Ann.} {\bf 346} (2010),  571--627.
%



\bibitem[Mu]{Mu} F.Murat,  Soluciones renormalizadas de EDP elipticas no lineales,
Preprint 93023, Laboratoire d'Analyse Num\'erique de
l'Universit\'e Paris VI (1993).


%
\bibitem[Scha]{Schauder} J.Schauder, 
Sur les \'equations lin\'eaires du type elliptique a coefficients continus
\emph{C. R. Acad. Sci. Paris} {\bf 199} (1934), 1366--1368.


\bibitem[Schw]{sebastian} S.Schwarzacher, Regularity for Degenerate Elliptic and Parabolic Systems, {\emph PhD thesis}, LMU Munich, Germany,  (2013).


\bibitem[Si,J.]{SimonJ} J.Simon
R\'egularit\'e de solutions de problemes nonlin\'eaires,
\emph{C. R. Acad. Sci. Paris S\'er. A-B} {\bf 282} (1976), A1351--A1354.

 




\bibitem[Uh]{Uhl} K.Uhlenbeck, Regularity for a class of non-linear elliptic
systems, \emph{Acta Math.} {\bf 138} (1977), 219--240.

\end{thebibliography}
\end{document}